\documentclass[twoside,11pt,leqno]{article}

\usepackage{amsfonts}

\newtheorem{theorem}{Theorem}[section]
\newtheorem{lemma}[theorem]{Lemma}

\newtheorem{statement}[theorem]{Statement}

\long\gdef\boxit#1{\begingroup\vbox{\hrule\hbox{\vrule\kern3pt
      \vbox{\kern3pt#1\kern3pt}\kern3pt\vrule}\hrule}\endgroup}

\def\qed{ \ \vrule width.2cm height.2cm depth0cm\smallskip}
\def\qedt{ \ \vrule width.2cm height.2cm depth0cm}
\def\Keywords{\bigskip\par {\sl Keywords\/}:\enspace}

\newenvironment{proof}{\noindent {\bf Proof.\/}}{$\qed$\vskip 0.1in}
\newcommand{\Xcomment}[1]{}

\newcommand{\sumn}{\sum\nolimits}

\newcommand{\refeq}[1]{(\ref{eq:#1})}

\newenvironment{myitem}{\refstepcounter{equation}\begin{enumerate}%
\item[(\thesection.\arabic{equation})]$\quad$}{\end{enumerate}}

\newenvironment{myitem1}{\refstepcounter{equation}\begin{enumerate}%
\item[(\thesection.\arabic{equation})]\begin{itemize}}%
{\end{itemize}\end{enumerate}}

\makeatletter
\@addtoreset{equation}{section}

\renewcommand{\section}{\@startsection{section}{1}{0pt}%
{-3.5ex plus -1ex minus -.2ex}{2.3ex plus .2ex}%
{\normalfont\Large}}

\renewcommand{\subsection}{\@startsection{subsection}{2}{0pt}%
{-3.25ex plus -1ex minus -.2ex}{1.5ex plus .2ex}%
{\normalfont\large\bf}}

\makeatother

\def\Rset{{\mathbb R}}
\def\Zset{{\mathbb Z}}
\def\Cset{{\mathbb C}}

\def\Cscr{{\cal C}}
\def\Dscr{{\cal D}}

\def\Fscr{{\cal F}}

\def\Kscr{{\cal K}}
\def\Lscr{{\cal L}}
\def\Pscr{{\cal P}}
\def\Rscr{{\cal R}}

\def\Tscr{{\cal T}}

\def\diver{\mbox{div}}
\def\supp{\mbox{supp}}

\def\tilde{\widetilde}
\def\bar{\overline}

\def\BC{{\cal B}}

\def\Pmin{\Pscr^{\min}}
\def\Kmin{\Kscr^{\min}}

\begin{document}
\begin{titlepage}

\title{\large\bf Concave Cocirculations in a Triangular Grid}

\def\thepage {} 

\author{{\large Alexander V. Karzanov}
\\
\ \\
{\sl Institute for System Analysis}\\
{\sl 9, Prospect 60 Let Oktyabrya, 117312 Moscow, Russia}\\
{{\sl E-mail}: {\it sasha@cs.isa.ac.ru}}
}

\date{March 2003}

\maketitle

\begin{abstract}
Let $G=(V(G),E(G))$ be a planar digraph embedded in the plane in
which all inner faces are equilateral triangles (with three edges
in each), and let the union
$\Rscr$ of these faces forms a convex polygon.
The question is: given a function $\sigma$ on the boundary edges of
$G$, does there exist a concave function $f$ on $\Rscr$ which is
affinely linear within each bounded face and satisfies
$f(v)-f(u)=\sigma(e)$ for each boundary edge $e=(u,v)$?

The functions $\sigma$ admitting such an $f$ form a
polyhedral cone $C$, and when the region $\Rscr$ is a triangle,
$C$ turns out to be exactly the cone of boundary data of honeycombs.
Studing honeycombs in connection with a problem on spectra of triples
of zero-sum Hermitian matrices, Knutson, Tao, and Woodward~\cite{KTW}
showed that $C$ is described by linear inequalities of Horn's type
with respect to so-called {\em puzzles}, along with obvious linear
constraints.

The purpose of this paper is to give an alternative proof of that
result, working in terms of discrete concave finctions, rather than
honeycombs, and using only linear programming and combinatorial tools.
Moreover, we extend the result to an arbitrary convex polygon $\Rscr$.
   \end{abstract}

\medskip
\Keywords Discrete convex function, Cocirculation, Honeycomb, Planar
graph, Flow

\bigskip

\end{titlepage}

\baselineskip 15pt

\section{Introduction} \label{sec:intr}

Let $\xi_1,\xi_2,\xi_3$ be three affinely independent vectors in the
euclidean plane $\Rset^2$ whose sum is the zero vector. The triangular
lattice generated by $\xi_1,\xi_2,\xi_3$ is associated with the infinite
planar directed graph $\Lscr$ whose vertices are integer combinations
of these vectors and whose edges are the ordered pairs $(u,v)$ of
vertices such that $v-u\in\{\xi_1,\xi_2,\xi_3\}$. An edge $(u,v)$ is
identified with the straight-line segment between $u,v$ oriented from
$u$ to $v$.

Consider a convex region $\Rscr$ in the plane formed by the union of
a nonempty finite set of faces ({\em little triangles}) of $\Lscr$;
it is a polygon with 3 to 6 sides.
We refer to the subgraph $G=(V(G),E(G))$ of $\Lscr$
consisting of the vertices and edges occurring in $\Rscr$ as a
{\em convex (triangular) grid}. The vertices and edges in the boundary
$b(G)$ of $G$ are called {\em outer} and their sets are denoted by
$V_0(G)$ and $E_0(G)$, respectively.
A real-valued function $f$ on the vertices of $G$ is called {\em
discrete concave (convex)} if it is the restriction to $V(G)$ of a
concave (resp. convex) function on $\Rscr$ which is affinely linear
within each bounded face of $G$. In this paper we prefer to deal with
discrete concave functions; the corresponding results for discrete
convex functions follow by symmetry.

We are interested in the functions on the set of
outer vertices that can be extended to discrete concave functions on
all vertices of $G$. Instead, one can consider the corresponding
functions on edges. More precisely, a function $h:E(G)\to\Rset$ is said to
be a {\em cocirculation} if it can be expressed via some function
$f:V(G)\to\Rset$ as $h(e):=f(v)-f(u)$ for each edge $e=(u,v)$. Such an $h$
determines $f$ up to a constant, and we refer to $h$ as a {\em concave
cocirculation} if $f$ is discrete concave. We consider the problem:
  \begin{myitem}
given a function $\sigma:E_0(G)\to\Rset$, decide whether $\sigma$
is extendable to a concave cocirculation in $G$.
  \label{eq:problem}
  \end{myitem}

Let $\BC(G)$ be the set of $\sigma$'s admitting such an extension.
Obvious necessary (but far to be sufficient) conditions on
$\sigma$ to belong to $\BC(G)$ are that the sum of values of $\sigma$,
taken with signs $+$ or $-$ depending on the direction of an edge in
the boundary circuit, amounts to zero and that $\sigma$ is monotone
nonincreasing along each side-path of $b(G)$.

When the region $\Rscr$ spanned by $G$ is a triangle, $G$ is called
a {\em 3-side (triangular) grid}. Then the boundary of $G$ is
the concatenation of three paths $B_1,B_2,B_3$ forming the sides
of $\Rscr$, where the edges of $B_i$ are parallel to $\xi_i$.
We say that $G$ is of {\em size} $n$ if $|B_i|=n$, where $|P|$
denotes the number of edges of a path $P$. In this special case
$\BC(G)$ coincides with the so-called Berenstein-Zelevinsky cone,
which also arises in other interesting models.
More precisely, relying on some earlier results, Knutson and
Tao~\cite{KT} showed that for a triple of monotone nonincreasing
$n$-tuples $(\lambda,\mu,\nu)\in(\Rset^n)^3$, the following properties
are equivalent:
\begin{itemize}
\item[(P1)] $\lambda,\mu,\nu$ are the spectra of three $n\times n$
Hermitian matrices whose sum is the zero matrix;
  \item[(P2)] there exists a {\em honeycomb} of size $n$ in which
the three tuples of semiinfinite edges have the constant coordinates
$\lambda,\mu,\nu$;
  \item[(P3)] let $G$ be the 3-side gride of size $n$ and let the
function $\sigma$ on $E_0(G)$ take the value $\lambda_j$ (resp.
$\mu_j,\nu_j$) on $j$th edge of the path $B_1$ (resp. $B_2,B_3$); then
$\sigma\in\BC(G)$.
  \end{itemize}

For a survey on those earlier results, see~\cite{Ful,KT}, and for
the definition of honeycomb, see~\cite{KT}.
Note that while the equivalence of (P2) and (P3) is rather transparent
(they are related via Fenchel's duality), the equivalence of
these to (P1) is quite sophisticated. In an earlier period of study of
property (P1) Horn~\cite{Ho} recursively constructed a finite list
of nontrivial necessary conditions on $\lambda,\nu,\mu$ to satisfy
this property and conjectured the sufficiency of this list (which,
in particular, implies that these $(\lambda,\mu,\nu)$'s constitute a
polyhedral cone). Horn's conditions are viewed as linear inequalities
of the form
  \begin{equation} \label{eq:horn}
\lambda(I)+\mu(J)+\nu(K)\ge 0
  \end{equation}
for certain subsets $I,J,K$ of $\{1,\ldots,n\}$ with $|I|=|J|=|K|$,
letting $\alpha(S):=\sum(\alpha_i: i\in S)$ for
$\alpha=(\alpha_1,\ldots,\alpha_n)\in\Rset^n$ and
$S\subseteq\{1,\ldots,n\}$. Subsequent efforts of several authors
(where the breakthrough contribution was due to Klyachko~\cite{Kl})
have resulted in a proof of Horn's conjecture; the obtained result is
referred in~\cite{KTW} as the ``H-R/T/K theorem'', abbreviating the
names of Helmke, Rosenthal, Totaro, and Klyachko. Recently Knutson,
Tao and Woodward~\cite{KTW} established a combinatorial existence
criterion for honeycombs, obtaining another proof of that theorem via
the above-mentioned relation to (P1). According to their criterion,
each Horn's triple $\{I,J,K\}$ is induced by a {\em puzzle}, a certain
subdivision of a 3-side grid into little triangles and little rhombi
endowed with a certain 0,1 weighting on the sides of these pieces.
An alternative method of deducing the characterization of the
$(\lambda,\mu,\nu)$'s in (P1) is suggested by Danilov and
Koshevoy~\cite{DK2}.

The purpose of this paper is to give a direct proof of the puzzle
criterion of the solvability of problem~\refeq{problem}, without
using relationships to honeycombs.
More precisely, we extend the notion of puzzle in a natural way to
an arbitrary convex grid $G$ and show that $\sigma:E_0(G)\to\Rset$
is extendable to a concave cocirculation if and only if it obeys the
linear inequalities of Horn's type determined by puzzles and the
above-mentioned obvious linear constraints. The proof uses only a
linear programming approach and combinatorial techniques.

This paper is organized as follows. Section~\ref{sec:prelim} contains
basic definitions and facts and states problem~\refeq{problem} as a
linear program.
In Section~\ref{sec:theorem} we explain the notion of puzzle
for a convex grid (using a somewhat different, but equivalent,
definition) and formulate the main theorem giving the puzzle criterion
of the solvability of the problem (Theorem~\ref{tm:puzzle}).
The proof of the theorem is given in
Section~\ref{sec:proof}, based on a weaker, linear programming,
criterion of the solvability of~\refeq{problem} established in
Sections~\ref{sec:lp} and on a representation of some dual variables
as flows in a certain dual graph, explained in Section~\ref{sec:flow}.
The concluding Section~\ref{sec:concl} considers a slightly more
general problem. It also discusses a generalization to convex
grids of results in~\cite{KTW} on the puzzles determining facets of
the cone $\BC(G)$ for a 3-side grid $G$, which characterize such
puzzles combinatorially and in terms of rigidity.

Related and other aspects of discrete convex (concave) functions on
triangular grids in the plane and some applications in analysis and
algebra are studied in~\cite{DK1}. For a survey on
discrete convex functions on lattices of higher dimensions and their
applications in combinatorics (in particular, in the theory of
matroids and submodular functions), see, e.g.,~\cite{Mur}.

\section{Preliminaries}
\label{sec:prelim}

\Xcomment{Let us call the functions in $\BC(G)$ {\em feasible}.}

We start with terminology, notation and conventions.
Edges, faces, subgraphs, paths, circuits and other
relevant objects in a convex grid $G$ or another graph in
question are usually identified with their {\em closed} images in the
plane. By a path (circuit) we usually mean a simple directed path
(circuit) $P=(v_0,e_1,v_1,\ldots,e_k,v_k)$, where $e_i=(v_{i-1},v_i)$;
it may be abbreviately denoted as $(e_1,e_2,\ldots,e_k)$ (via edges).
A path $P$ with beginning vertex $u$ and end vertex $v$ is called a
$u-v$~{\em path}; $P$ is called {\em degenerate} if $u=v$ (i.e.,
$|P|=0$).
When $P$ forms a straight-line segment in the plane, $P$ is
called a {\em straight path}, or a {\em line} of the graph.
A $k$-{\em circuit} is a circuit with $k$ edges.

Since problem~\refeq{problem} does not depend, in essense, on the
choice of vectors $\xi_1,\xi_2,\xi_3$ (subject to the condition in the
Introduction), we fix these for convenience as
$\xi_1=(1,0)$, $\xi_2=(-1,\sqrt{3})/2$ and $\xi_3=(-1,-\sqrt{3})/2$.
Then the little triangles of $G$ are equilateral triangles of
{\em size} 1. Note that the boundary of any triangle in
$G$ (formed by the union of some faces) is a circuit directed
clockwise or counterclockwise around the triangle. The little triangle
surrounded by a 3-circuit $C$ is denoted by $\Delta_C$.
We say that a triangle is {\em normal} if its boundary circuit is
directed counterclockwise, and {\em turned-over} otherwise.
 \begin{center}
  \unitlength=1mm
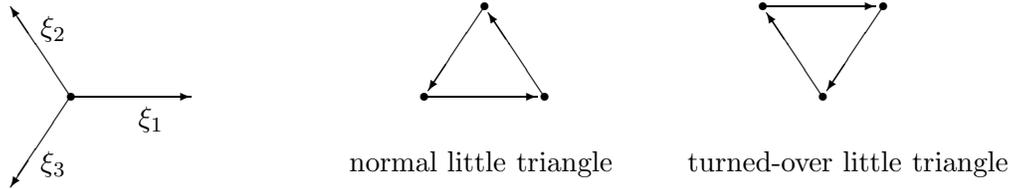
  \begin{picture}(130,24)
\put(8,12){\circle*{1.0}}
\put(8,12){\vector(1,0){16}}
\put(8,12){\vector(-2,-3){8}}
\put(8,12){\vector(-2,3){8}}
\put(17,8){$\xi_1$}
\put(4,20){$\xi_2$}
\put(4,2){$\xi_3$}
\put(55,12){\circle*{1.0}}
\put(63,24){\circle*{1.0}}
\put(71,12){\circle*{1.0}}
\put(55,12){\vector(1,0){15}}
\put(63,24){\vector(-2,-3){7.5}}
\put(71,12){\vector(-2,3){7.5}}
\put(45,2){normal little triangle}
\put(108,12){\circle*{1.0}}
\put(100,24){\circle*{1.0}}
\put(116,24){\circle*{1.0}}
\put(100,24){\vector(1,0){15}}
\put(116,24){\vector(-2,-3){7.5}}
\put(108,12){\vector(-2,3){7.5}}
\put(90,2){turned-over little triangle}
  \end{picture}
 \end{center}

We denote the sets of outer edges directed counterclockwise and
clockwise around $\Rscr$ by $E^+_0(G)$ and $E^-_0(G)$,
respectively. A maximal straight path in $b(G)$, or a
{\em side-path} of $G$, whose edges are parallel to
$\xi_i$ and belong to $E^+_0(G)$ (resp. $E^-_0(G)$) is denoted
by $B^+_i$ (resp. $B^-_i$). One may assume that if $G$ is a 3-side
grid, then the boundary of $G$ is formed by $B^+_1,B^+_2,B^+_3$.

For a function $h$ on $E(G)$, its restriction to the set of outer
edges is called the {\em border} of $h$.

Next we explain how to write problem~\refeq{problem} as a linear
program. Obviously, a function $f:V(G)\to\Rset$ is discrete concave if
and only if
 \begin{equation} \label{eq:concave}
f(u)+f(u')\le f(v)+f(v')
  \end{equation}
holds for each {\em little rhombus} (the union of two little
triangles sharing a common edge) $\rho$, where $u,u'$ are the
{\em acute} vertices and $v,v'$ are the {\em obtuse} vertices of
$\rho$:
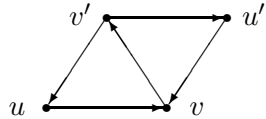
 \begin{center}
  \unitlength=1mm
  \begin{picture}(35,13)
\put(5,1){\circle*{1.0}}
\put(21,1){\circle*{1.0}}
\put(13,13){\circle*{1.0}}
\put(29,13){\circle*{1.0}}
\put(0,0){$u$}
\put(24,0){$v$}
\put(8,12){$v'$}
\put(31,12){$u'$}
\put(5,1){\vector(1,0){15.5}}
\put(13,13){\vector(1,0){15.5}}
\put(13,13){\vector(-2,-3){7.7}}
\put(29,13){\vector(-2,-3){7.7}}
\put(21,1){\vector(-2,3){7.7}}
  \end{picture}
 \end{center}

Clearly $h\in\Rset^{E(G)}$ is a cocirculation if and only if
the sum of its values on each 3-circuit is zero. (Hereinafter we do
not distinguish between real-valued functions on a finite set $E$
and vectors in the Euclidean space $\Rset^E$ whose coordinates are
indexed by the elements of $E$.) Linear
constraints reflecting the property of a cocirculation $h$ to be
concave are derived from~\refeq{concave}.
Let us say that an ordered pair $\tau=(e,e')$ of non-adjacent edges of
$G$ is a {\em tandem} if they occur as opposite sides
of a little rhombus $\rho$ and the head of $e$ is an obtuse
vertex of $\rho$ (while the other obtuse vertex of $\rho$ is
the tail of $e'$). We distinguish between two sorts of tandems by
specifying $\tau$ as a {\em normal} tandem if the little triangle in
$\rho$ containing $e$ is normal, and a {\em turned-over} tandem
otherwise. Note that each little rhombus $\rho$ involves two tandems
one of which is normal and the other is turned-over. The picture
illustrates the case when $e,e'$ are parallel to $\xi_1$.
 \begin{center}
  \unitlength=1mm
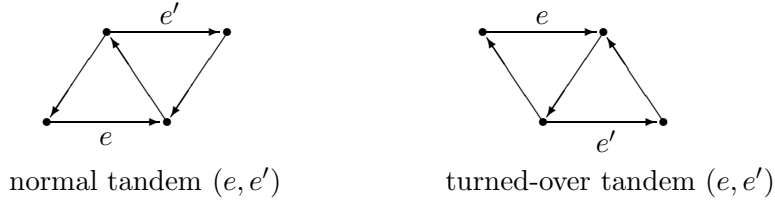
  \begin{picture}(95,26)
\put(5,10){\circle*{1.0}}
\put(21,10){\circle*{1.0}}
\put(13,22){\circle*{1.0}}
\put(29,22){\circle*{1.0}}
\put(12,7){$e$}
\put(20,23){$e'$}
\put(0,1){normal tandem $(e,e')$}
\put(5,10){\vector(1,0){15}}
\put(13,22){\vector(1,0){15}}
\put(13,22){\vector(-2,-3){7.5}}
\put(29,22){\vector(-2,-3){7.5}}
\put(21,10){\vector(-2,3){7.5}}
\put(87,10){\circle*{1.0}}
\put(71,10){\circle*{1.0}}
\put(63,22){\circle*{1.0}}
\put(79,22){\circle*{1.0}}
\put(78,6){$e'$}
\put(70,23){$e$}
\put(58,1){turned-over tandem $(e,e')$}
\put(71,10){\vector(1,0){15}}
\put(63,22){\vector(1,0){15}}
\put(79,22){\vector(-2,-3){7.5}}
\put(71,10){\vector(-2,3){7.5}}
\put(87,10){\vector(-2,3){7.5}}
\Xcomment{            right fig.
\put(110,10){\circle*{1.0}}
\put(126,10){\circle*{1.0}}
\put(118,22){\circle*{1.0}}
\put(117,7){$e$}
\put(109,15){$e''$}
\put(123,15){$e'$}
\put(100,1){3-circuit $(e,e',e'')$}
\put(110,10){\vector(1,0){15}}
\put(118,22){\vector(-2,-3){7.5}}
\put(126,10){\vector(-2,3){7.5}}
}
  \end{picture}
 \end{center}

For the cocirculation $h$ generated by a function $f$ on the vertices,
\refeq{concave} is just equivalent to the condition $h(e)\ge h(e')$
on the normal tandem $(e,e')$ in the little rhombus $\rho$. Thus,
given $\sigma\in \Rset^{E_0(G)}$, a concave cocirculation with the
border $\sigma$ is a solution $h\in\Rset^{E(G)}$ of the system:
\begin{equation} \label{eq:lp1}
     h(e)+h(e')+h(e'')=0, \qquad C=(e,e',e'')\in\Cscr(G),
\end{equation}
\vspace{-0.5cm}
\begin{equation} \label{eq:lp2}
    h(e')-h(e)\le 0, \qquad \tau=(e,e')\in\Tscr(G),
\end{equation}
\vspace{-0.5cm}
\begin{equation} \label{eq:lp3}
    h(e)=\sigma(e), \qquad e\in E_0(G),
\end{equation}
where $\Cscr(G)$ is the set of 3-circuits (considered
up to cyclically shifting), and $\Tscr(G)$ the set of normal tandems
in $G$. When this system has a solution, we call $\sigma$ {\em
feasible}.

\smallskip
As mentioned in the Introduction, there are two elementary conditions
on $\sigma$ to be feasible. The first one (necessary for the border of
any cocirculation) is the {\em zero-sum condition}:
 \begin{equation}   \label{eq:zero-sum}
\sigma(E^+_0(G))-\sigma(E^-_0(G))=0,
 \end{equation}
The second one is the {\em monotone condition}:
 \begin{myitem}
$\sigma(e)\ge\sigma(e')$ for each straight path $(e,e')$ of $b(G)$.
  \label{eq:monot}
  \end{myitem}
 \Xcomment{
This follows by considering the edge $q\in E(G)$ such that both $(e,q)$
and $(q,e')$ are tandems (such a $q$ exists and unique). Then $h(e)\ge
h(q)\ge h(e')$, by~\refeq{lp2}.  }

Since the set of concave cocirculations on $G$ is described by a
finite number of linear constraints, the set $\BC(G)$ of feasible functions
$\sigma$ (the borders of concave cocirculations in $G$) forms a
convex polyhedral cone in $\Rset^{E_0(G)}$.
To compute the dimension of this cone is
easy (cf.~\cite{KTW}).
  \begin{statement}  \label{st:dimB}
$\dim(\BC(G))=|E_0(G)|-1$.
  \end{statement}
\begin{proof}
In view of~\refeq{zero-sum}, $\dim(\BC(G))\le|E_0(G)|-1=:r$.
To show the reverse inequality, we first construct a concave
cocirculation $h$ for which all tandem inequalities in~\refeq{lp2} are
strict.

Take a maximal straight $u-v$~path $P$ of $G$ not
contained in $b(G)$. Let $Z$ be the set of edges of $G$ that lie in the
region on the right from $P$ (when moving from $u$ to $v$) and are
not parallel to $P$. Define $h_P(e)$ to be 1 if $e\in Z$ and $e$
points toward $P$, $-1$ for the other edges $e$ in $Z$, and 0 for the
remaining edges of $G$. One can check that $h_P$ is a concave
cocirculation and that $h(e)>h(e')$ for each tandem $(e,e')$ where $e$
and $e'$ are separated by $P$. The sum of $h_P$'s over all such paths
$P$ gives the desired concave cocirculation $h$. Let $\sigma$ be the
border of $h$.

Now for each outer vertex $v$ and each edge $e$, define $h_v(e)$ to be
1 if $v$ is the head of $e$, $-1$ if $v$ is the tail of $e$, and 0
otherwise. Then $h_v$ is a cocirculation; moreover,
$h+\frac12 h_v$ is a concave cocirculation. Let
$\sigma_v$ be the border of $h_v$. Clearly $r$ borders
among these $\sigma_v$ are linearly independent. This implies that
$r$ borders $\sigma+\frac12 \sigma_v$ of the concave cocirculations
$h+\frac12 h_v$ are linearly independent.
  \end{proof}

\section{Theorem}
\label{sec:theorem}

Linear programming (Farkas lemma) suggests a standard way to obtain a
solvability criterion for system~\refeq{lp1}--\refeq{lp3}. Our aim,
however, is to obtain a sharper, combinatorial, characterization for
the borders of concave cocirculations on $G$.

First of all we construct a certain dual digraph $H$. For each edge
$e\in E(G)$, take the median point $v_e$ on $e$, making it a
vertex of $H$. For each normal tandem $\tau=(e,e')$,
form (straight-line) edge $a_\tau$ from $v_e$ to $v_{e'}$, making it
an edge of $H$.
Note that when $e,e'$ are parallel to $\xi_i$, the edge $a_\tau$ is
{\em anti-parallel} to $\xi_{i-1}$, in the sense that $a_\tau$ is a
parallel translate of the opposite vector $-\xi_{i-1}$. (Hereinafter
the corresponding indices are taken modulo 3.)
The resulting graph $H$ is the union of three disjoint digraphs
$H_1,H_2,H_3$, where $H_i$ is induced by the introduced edges
connecting points on edges of $G$ parallel to $\xi_i$. The three types
of edges of $H$ are drawn in bold in the picture.
 \begin{center}
  \unitlength=1mm
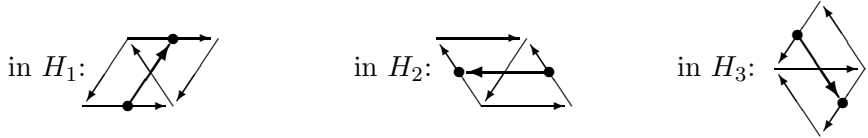
  \begin{picture}(135,22)
\put(10,5){\begin{picture}(18,9)%
  \put(0,0){\vector(1,0){11}}
  \put(6,9){\vector(1,0){11}}
  \put(12,0){\vector(-2,3){5.5}}
  \put(6,9){\vector(-2,-3){5.5}}
  \put(18,9){\vector(-2,-3){5.5}}
  \end{picture}}
\put(0,9){in $H_1$:}
\put(57,5){\begin{picture}(18,9)%
  \put(6,0){\vector(1,0){11}}
  \put(0,9){\vector(1,0){11}}
  \put(6,0){\vector(-2,3){5.5}}
  \put(18,0){\vector(-2,3){5.5}}
  \put(12,9){\vector(-2,-3){5.5}}
  \end{picture}}
\put(46,9){in $H_2$:}
\put(102,1){\begin{picture}(12,18)%
  \put(0,9){\vector(1,0){11}}
  \put(6,0){\vector(-2,3){5.5}}
  \put(12,9){\vector(-2,3){5.5}}
  \put(6,18){\vector(-2,-3){5.5}}
  \put(12,9){\vector(-2,-3){5.5}}
  \end{picture}}
\put(89,9){in $H_3$:}
\put(16,5){\circle*{1.5}}
\put(22,14){\circle*{1.5}}
\put(72,9.5){\circle*{1.5}}
\put(60,9.5){\circle*{1.5}}
\put(105,14.5){\circle*{1.5}}
\put(111,5.5){\circle*{1.5}}
\thicklines
\put(16,5){\vector(2,3){5.5}}
\put(72,9.5){\vector(-1,0){11}}
\put(105,14.5){\vector(2,-3){5.5}}
  \end{picture}
 \end{center}

So the maximal paths in $H_i$ are straight, pairwise disjoint and
anti-parallel to $\xi_{i-1}$. If a path $P$
of $H$ begins at $v_e$ and ends at $v_{e'}$, we say that $P$ {\em
leaves} the edge $e$ and {\em enters} the edge $e'$ (of $G$),
admitting the case of degenerate $P$. We also say that $P$ leaves
(enters) a little triangle $\Delta$ if $e\subset\Delta$ (resp.
$e'\subset\Delta$).

\medskip
{\bf Definition.} A {\em puzzle} is a pair $\Pi=(\Fscr,\Pscr)$
consisting of a set $\Fscr$ of little triangles of $G$ and a set
$\Pscr$ of paths of $H$ such that:
  \begin{myitem1}
\item[(i)] the interiors of triangles in $\Fscr$ and paths in
$\Pscr$ are pairwise disjoint;
  \item[(ii)] for each edge $e$ of each normal (resp. turned-over)
triangle in $\Fscr$, there is precisely one path in $\Pscr$ entering
(resp. leaving) $e$;
  \item[(iii)]
for each path in $\Pscr$ leaving edge $e$ and entering edge $e'$,
either $e$ belongs to a turned-over triangle in $\Fscr$ or
$e\in E^+_0(G)$, and similarly, either $e'$ belongs to a normal
triangle in $\Fscr$ or $e\in E^-_0(G)$.
  \label{eq:def_gal}
  \end{myitem1}

\smallskip
(Degenerate paths $P=v_e$ in $\Pscr$ are admitted. When $e$ is an
inner edge, such a $P$ serves to ``connect''the pair of triangles in
$\Fscr$ sharing the edge $e$. When $e$ is an outer edge, $P$
``connects'' this edge with the triangle in $\Fscr$ containing $e$.)
The {\em boundary} $b(\Pi)$ of $\Pi$ is defined to be set of
outer edges $e$ for which there is a path in $\Pscr$ leaving
or entering $e$. The subsets of edges of $E^+_0(G)$ and $E^-_0(G)$
occurring in $b(\Pi)$ are denoted by $b^+(\Pi)$ and $b^-(\Pi)$),
respectively.

The {\em puzzle criterion} of
the solvability of~\refeq{lp1}--\refeq{lp3} is the following.

\begin{theorem} \label{tm:puzzle}
Let $G$ be a convex triangular grid, and let $\sigma:E_0(G)\to\Rset$
satisfy~\refeq{zero-sum} and~\refeq{monot}. Then a concave
cocirculation $h$ in $G$ with $h(e)=\sigma(e)$ for all $e\in E_0(G)$
exists if and only if
  \begin{equation}  \label{eq:gall-ineq}
\sigma(b^+(\Pi))-\sigma(b^-(\Pi)) \ge 0
  \end{equation}
holds for each puzzle $\Pi$.
  \end{theorem}

Thus, the cone $\BC(G)$ is described by the
{\em puzzle inequalities}~\refeq{gall-ineq} and the linear
constraints~\refeq{zero-sum} and~\refeq{monot}.

\medskip
\noindent {\bf Remark.}
A puzzle in a 3-side grid $G$ introduced in Knutson et
al.~\cite{KTW} is defined to be a diagram $D$ consisting of a
subdivision of the big triangle $\Rscr$ into little triangles and
little rhombi of $G$, and of a 0,1 labelling of the edges of $G$ that
are sides of these pieces,
satisfying the following conditions: (a) the three sides of each
little triangle are labelled either 1,1,1 or 0,0,0, and (b) the sides
of each little rhombus $\rho$ are labelled 0,1,0,1, in this order
clockwise of an acute vertex of $\rho$ (the
{\em triangle-rhombus-label} form). The boundary $b(D)$ of $D$ is
defined to be the set of outer edges labelled 1.
There is a natural one-to-one correspondence between the puzzles $D$
of~\cite{KTW} and those $\Pi=(\Fscr,\Pscr)$ in the above definition
(the {\em triangle-path} form)
and this correspondence preserves the puzzle boundary: $b(D)=b(\Pi)$.
(In this correspondence, $\Fscr$ is set of little triangles labelled
1,1,1, and the edges of $H$ used in the paths of $\Pscr$ are those
connecting the sides labelled 1 in the rhombi of $D$.)
The triangle-path form of puzzle is more convenient for us to handle
in the proof of Theorem~\ref{tm:puzzle}, which is based on certain
path and flow constructions.

\medskip
To illustrate the theorem, consider a 3-side grid of size $n$ and
a puzzle having one triangle $\Delta$ and three paths $P_1,P_2,P_3$,
each $P_i$ connecting $\Delta$ with the side-path
$B^+_i=(b_i^1,\ldots,b_i^n)$.
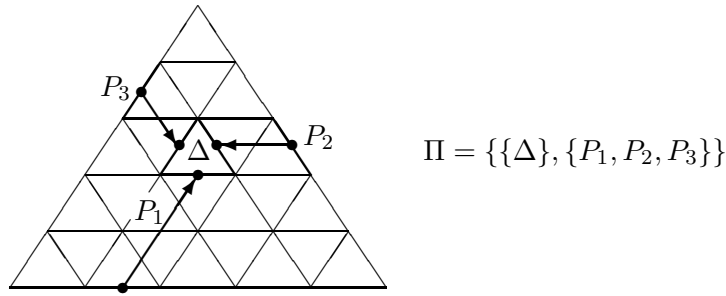
 \begin{center}
  \unitlength=1mm
  \begin{picture}(90,38)
\put(0,0){\line(1,0){50}}
\put(5,7.5){\line(1,0){40}}
\put(10,15){\line(1,0){30}}
\put(15,22.5){\line(1,0){20}}
\put(20,30){\line(1,0){10}}
\put(0,0){\line(2,3){25}}
\put(10,0){\line(2,3){6}}   
\put(30,30){\line(-2,-3){11.5}}   
\put(20,0){\line(2,3){15}}
\put(30,0){\line(2,3){10}}
\put(40,0){\line(2,3){5}}
\put(10,0){\line(-2,3){5}}
\put(20,0){\line(-2,3){10}}
\put(30,0){\line(-2,3){15}}
\put(40,0){\line(-2,3){20}}
\put(50,0){\line(-2,3){25}}
\thicklines
\put(10,0){\line(1,0){10}}
\put(15,22.5){\line(2,3){5}}
\put(40,15){\line(-2,3){5}}
\put(20,15){\line(1,0){10}}
\put(20,15){\line(2,3){5}}
\put(30,15){\line(-2,3){5}}
\put(15,0){\vector(2,3){9.7}}
\put(37.5,19){\vector(-1,0){9.5}}
\put(17.5,26){\vector(2,-3){4.7}}
\thinlines
\put(15,0){\circle*{1.5}}
\put(25,15){\circle*{1.5}}
\put(17.5,26){\circle*{1.5}}
\put(22.5,19){\circle*{1.5}}
\put(27.5,19){\circle*{1.5}}
\put(37.5,19){\circle*{1.5}}
\put(16.5,9){$P_1$}
\put(39,19){$P_2$}
\put(12,25.5){$P_3$}
\put(23.5,17){$\Delta$}
\put(55,17){$\Pi=\{\{\Delta\},\{P_1,P_2,P_3\}\}$}
  \end{picture}
 \end{center}
\noindent
Let $P_i$ leave edge $b_i^{r(i)}\in B_i$ and enter edge
$e_i\subset \Delta$. Summing up the inequalities in~\refeq{lp2}
for the normal tandems induced by the edges of $P_i$, we have
$\sigma(b_i^{r(i)})=h(b_i^{r(i)})\ge h(e_i)$. This together
with~\refeq{lp1} for the 3-circuit $(e_1,e_2,e_3)$ implies
$\sum(\sigma(b_i^{r(i)}):i=1,2,3)\ge 0$. Also $r(1)+r(2)+r(3)=n+2$.
Thus, any feasible $\sigma=(\lambda,\mu,\nu)\in(\Rset^n)^3$ must obey
  $$
\lambda_i+\mu_j+\nu_k\ge 0
  $$
for any choice of $i,j,k$ with $i+j+k=n+2$.
This is the simplest sort of Horn's inequality~\refeq{horn}.

\smallskip
One can associate with a puzzle $\Pi=(\Fscr,\Pscr)$ graph
$\Gamma_\Pi$ whose vertices are the triangles in $\Fscr$ and the edges
in $b(\Pi)$ and where the vertices $u,v$ are connected by an edge
if and only if there is a path in $\Pscr$ leaving one and entering the
other of $u,v$ (such a $\Gamma_\Pi$ describes the topological type of
$\Pi$, in a sense). It is not difficult to deduce that $\Gamma_\Pi$ is
determined by the list of cardinalities $|b(\Pi)\cap B|$, where $B$
ranges over the side-paths of $G$. In particular,
\begin{myitem}
the numbers $|\Fscr|$ and $|\Pscr|$ are determined by $b(\Pi)$.
 \label{eq:par_puzzle}
 \end{myitem}
(Instruction: shrink into a point each little triangle of $G$
that neither belongs to $\Fscr$ nor meets a path in $\Pscr$, and
simultaneously, for each nondegenerate $v_e-v_{e'}$ path in $\Pi$,
shrink into $e$ the parallelogram with opposite sides $e,e'$.
The resulting graph $G'$ is again a convex grid (possibly degenerate)
in which the little
triangles one-to-one correspond to those in $\Fscr$, and the edges to
the paths $\Pscr$; also the outer edges of $G'$ one-to-one correspond
to the edges in $b(\Pi)$ when $\Fscr\ne\emptyset$.
Moreover, $G'$ depends only on the above-mentioned cardinalities.)

\section{Linear Programming Approach}
             \label{sec:lp}

In what follows, speaking of a tandem, we always mean a normal
tandem in $G$.
Assign a variable $z(C)\in\Rset$ to each 3-circuit $C$ of $G$,
a variable $g(\tau)\in\Rset_+$ to each tandem $\tau$,
and a variable $d(e)\in\Rset$ to each outer edge $e$.
Then the linear system dual of~\refeq{lp1}--\refeq{lp3}
is viewed as
  \begin{equation}  \label{eq:dlp1}
\sum_{C\in\Cscr(G): e\in C}z(C) -\sum_{\tau=(e,e')\in\Tscr(G)}g(\tau)
 +\sum_{\tau=(e',e)\in\Tscr(G)}g(\tau)=0, \quad e\in E(G)-E_0(G),
  \end{equation}
\vspace{-0.4cm}
  \begin{equation}  \label{eq:dlp2}
\sum_{C\in\Cscr(G): e\in C}z(C) -\sum_{\tau=(e,e')\in\Tscr(G)}g(\tau)
 +\sum_{\tau=(e',e)\in\Tscr(G)}g(\tau) +d(e)=0, \quad e\in E_0(G).
  \end{equation}

Applying Farkas lemma to~\refeq{lp1}--\refeq{lp3},
we obtain the l.p. solvability criterion.
  \begin{statement}  \label{st:lin-crit}
Let $\sigma\in\Rset^{E_0(G)}$. A concave cocirculation $h$ with the
border $\sigma$ exists if and only if
 \begin{equation}  \label{eq:sigma-d}
\sigma\cdot d\ge 0
  \end{equation}
holds for any $z:\Cscr(G)\to\Rset$, $g:\Tscr(G)\to\Rset_+$ and
$d:E_0(G) \to\Rset$ satisfying~\refeq{dlp1} and~\refeq{dlp2}.
\qedt
   \end{statement}

Hereinafter for $a,b\in\Rset^E$, $a\cdot b$ denotes the inner product
$\sum(a(e)b(e):e\in E)$. We call a triple $K=(z,g,d)$
satisfying~\refeq{dlp1}--\refeq{dlp2} a {\em vector configuration},
or, briefly, a {\em v-configuration}, and regard $d$ as its border.

Statement~\ref{st:lin-crit} implies that the cone $\Dscr$ of borders
of v-configurations (which is convex) is anti-polar to the cone
$\BC(G)$ of borders of concave cocirculation in $G$, i.e.,
$\Dscr:=\{d\in\Rset^{E_0(G)}: \sigma\cdot d\ge 0\;
\forall \sigma\in\BC(G)\}$.
For an outer edge $e$, define $\theta(e):=1$ if $e\in E^+_0(G)$,
and $-1$ if $e\in E^-_0(G)$. Since the dimension of $\BC(G)$
is $|E_0(G)|-1$ (by Statement~\ref{st:dimB}) and $\BC(G)$
is contained in the hyperplane $\theta^\bot$ orthogonal to $\theta$
(by~\refeq{zero-sum}), the cone $\Dscr$ is full-dimensional and
contains the line $\Rset\theta$. So the facets of $\BC(G)$ one-to-one
correspond (by the orthogonality) to the 2-dimensional faces of
$\Dscr$, each being of the form $r_1 d+r_2\theta$
($r_1\in\Rset_+,r_2\in\Rset$) for a certain $d\in\Rset^{E_0(G)}$.

For a function (vector) $x$, let $\supp^+(x)$ and $\supp^-(x)$ denote
the {\em positive} part $\{e: x(e)>0\}$ and the {\em negative} part
$\{e: x(e)<0\}$ of the {\em support} $\supp(x)$ of $x$, respectively.
Since inequality~\refeq{sigma-d} is invariant in adding to $d$ any
multiple of $\theta$, it suffices to verify this inequality
only for the v-configurations $K=(z,g,d)$
satisfying:
  \begin{myitem}
(a) $\supp^+(d)\subseteq E^+_0(G)$ and $\supp^-(d)\subseteq
E^-_0(G)$, and (b) $\supp(d)\ne \emptyset,E_0(G)$.
  \label{eq:prop_conf}
  \end{myitem}
In what follows, we throughout assume that any v-configuration in
question satisfies (a). When (b) takes place too, we
call $K$ {\em proper}.

Let $\Sigma(G)$ be the set of $\sigma\in\Rset^{E_0(G)}$
satisfying~\refeq{zero-sum}--\refeq{monot}.
Then $\BC(G)\subseteq\Sigma(G)$.
A v-configuration $K=(z,g,d)$ is called {\em essential} if
$d$ separates $\Sigma(G)$, i.e., $\sigma\cdot d<0$ for some
$\sigma\in\Sigma(G)$.
Consider two v-configurations $K=(z,g,d)$ and $K'=(z',g',d')$
(subject to (a) in~\refeq{prop_conf}). $K$ and $K'$ are
called {\em equivalent} if their borders are proportional, i.e.,
$d=rd'$ for some $r>0$. We say that $K'$ {\em dominates} $K$ if at
least one of the following takes place:
  \begin{myitem1}
\item[(i)] $\sigma\in\Sigma(G)$ and $\sigma\cdot d<0$ imply
$\sigma\cdot d'<0$, and there exists $\sigma\in\Sigma(G)$ such
that $\sigma\cdot d\ge 0$ but $\sigma\cdot d'<0$; or
   \item[(ii)] $K$ is proper and not equivalent to $K'$,
and $K-rK'$ is a v-configuration (subject to (a)
in~\refeq{prop_conf}) for some $r>0$.
  \label{eq:dominat}
  \end{myitem1}

If $K$ is dominated by some $K'$, then $K$ is redundant and can be
excluded from consideration (as $d$ is \'a priori not
facet-determining for $\BC(G)$). This is obvious in case (i).
And in case (ii), the border $d'':=d-rd'$ of
the v-configuration $K'':=K-rK'$ is nonzero and
satisfies $(d'')^\bot \cap\BC(G)\supseteq d^\bot\cap\BC(G)$ and
$\supp(d'')\subseteq\supp(d)$. The former inclusion implies that if
$d^\bot$ contains a facet $F$ of $\BC(G)$, then $(d'')^\bot$ contains
$F$ as well. Then
$d''=r_1 d+r_2\theta$ for some $r_1>0$ and $r_2\in\Rset$, which
contradicts the latter inclusion since $\supp(d)\ne E_0(G)$ and
$K,K'$ are not equivalent.

Our method of proof of Theorem~\ref{tm:puzzle} consists in examining
an arbitrary essential
configuration $K$ and attempting to show that $K$ is
dominated unless it is equivalent to some ``puzzle configuration''.
Note that one can consider only rational-valued $z,g,d$
in~\refeq{dlp1}--\refeq{dlp2}. Moreover, by scaling, it suffices to
deal with {\em integer} v-configurations $(z,g,d)$.

For an outer edge $e$ of $G$, let $\chi^e$ denote the unit base vector
of $e$ in $\Rset^{E_0(G)}$ (i.e., $\chi^e(a)=1$ for $a=e$, and 0 otherwise).
We will use the following observation:
\begin{myitem}
if $K$ is an essential v-configuration with border $d$, $K'$ is a
v-configuration with border $d'$, and $d'=d-\chi^e+\chi^{e'}$, where
$e,e'$ are outer edges occurring in a straight path of $b(G)$ in this
order, then $K'$ dominates $K$.
  \label{eq:eep}
  \end{myitem}
To see this, let $d'':=\chi^e-\chi^{e'}$. Then $\sigma\cdot d''\ge 0$
for all $\sigma\in\Sigma(G)$, by~\refeq{monot}. This and
$d=d'+d''$ imply $\sigma\cdot d\ge 0$ for all $\sigma\in\Sigma(G)$
satisfying $\sigma\cdot d'\ge 0$.
Take $\sigma_1\in\BC(G)$ such that $\sigma_1(e)>\sigma_1(e')$
(existing by Statement~\ref{st:dimB}). Then $\sigma_1\cdot d'\ge 0$
and $\sigma_1\cdot d''>0$, implying $p:=\sigma_1\cdot d>0$. Take
$\sigma_2\in\Sigma(G)$ such that $q:=\sigma_2\cdot d<0$ (existing as
$K$ is essential).
Define $\sigma:=\sigma_2-\frac{q}{p}\sigma_1$.
We have $\sigma\cdot d=\sigma_2\cdot d-\frac{q}{p}\sigma_1\cdot
d=q-q=0$ and $\sigma\cdot d'=\sigma\cdot d-\sigma\cdot d''
=-\sigma\cdot d''=-\sigma_2\cdot d''+\frac{q}{p}\sigma_1\cdot
d''<0$, yielding~\refeq{dominat}(i).

\section{Flow Model} \label{sec:flow}

In the proof of Theorem~\ref{tm:puzzle} we will take advantage of a
representation of a v-configuration $K=(z,g,d)$ in a more
combinatorial form introduced in this section.

For a 3-circuit $C$, let us interprete $z(C)$ as the {\em weight}
of the little triangle $\Delta_C$. Similarly, $d(e)$ is the weight of
an outer edge $e$.
For each tandem $\tau=(e,e')$, set $g(a_\tau):=g(\tau)$,
interpreting it as the value of {\em flow} from $v_e$ to
$v_{e'}$ on the edge $a_\tau$ of the graph $H$ (introduced in
Section~\ref{sec:theorem}). The outer edges and little triangles
with nonzero weights are interpreted as ``sources'' or
``sinks'' of the flow. We say that an outer edge $e$ {\em emits}
$d(e)$ (units of) flow if $d(e)>0$, and {\em absorbs} $|d(e)|$ flow if
$d(e)<0$. Similarly, a little triangle $\Delta_C$ emits $z(C)$ flow
(through each of its three sides) if $z(C)>0$, and absorbs $|z(C)|$
flow if $z(C)<0$. Relations~\refeq{dlp1}--\refeq{dlp2} turn into the
flow balance condition
  \begin{equation}  \label{eq:flow}
\diver_g(v)+\sumn_{C\in\Cscr(G):v\in\Delta_C} z(C)
   +\sumn_{e\in E_0(G):v\in e} d(e)=0\qquad
          \mbox{for each}\;\; v\in V(H),
  \end{equation}
where
  $$
\diver_g(v):=\sumn_{u:(u,v)\in E(H)} g(u,v)-
               \sumn_{w:(v,w)\in E(H)} g(v,w).
  $$

Next, for a path $P$ in $H$, let $\chi^P\in \Rset^{E(H)}$ denote the
incidence vector of the set of edges of $P$.
Considering $g$ as a function on $E(H)$, applying
usual flow decomposition techniques and taking into
account~\refeq{flow}, one can find paths $P_1,\ldots,P_k$ in $H$
(possibly including degenerate paths) and positive real {\em weights}
$\alpha_1,\ldots,\alpha_k$ of these paths such that:
  \begin{myitem}
$g=\alpha_1\chi^{P_1}+\ldots+\alpha_k\chi^{P_k}$;
  \label{eq:path_form}
  \end{myitem}
  \begin{myitem}
for each edge $e$ of $G$, the sum of weights of emitting
elements containing $e$ is equal to the sum of weights of paths $P_i$
leaving $e$; similarly, the sum of absolute values of weights of
absorbing elements containing $e$ is equal to the sum of weights of
paths $P_i$ entering $e$.
  \label{eq:bal_paths}
  \end{myitem}
We call $(P_1,\ldots,P_k; \alpha_1,\ldots,\alpha_k)$
satisfying~\refeq{path_form}--\refeq{bal_paths} a {\em paths
decomposition} of $g$.

When $g$ is integer-valued, there is a decomposition
with all weights $\alpha_i$ integer
(an {\em integer paths decomposition}).
In this case we define a triple $\Kscr=(\Phi,\Pscr,\iota)$
{\em representing} $K$, in a sense, as follows.
Take $d(e)$ copies of each emitting outer edge $e$ and $z(C)$ copies of
each emitting triangle $\Delta_C$, forming family $\Phi^+$ of
(unweighted) emitting elements. Take $|d(e)|$ copies of each absorbing
outer arc $e$ and $|z(C)|$ copies of each absorbing
triangle $\Delta_C$, forming family $\Phi^-$ of absorbing
elements. Then $\Phi$ is the disjoint union of $\Phi^+$ and $\Phi^-$.
Take $\alpha_i$ copies of each path $P_i$, forming $\Pscr$. Assign
a map $\iota:\Pscr\to \Phi^+\times\Phi^-$ so as to satisfy the
following property:
 \begin{myitem}
if $P\in\Pscr$ and $\iota(P)=(\phi,\phi')$, then $P$ leaves $\phi$
and enters $\phi'$; moreover, for each $\phi\in\Phi^+$ (resp.
$\phi\in\Phi^-$) and each edge $e$ in $\phi$, there is exactly one
path $P\in\Pscr$ such that $\iota(P)=(\phi,\cdot)$ and $P$ leaves $e$
(resp. $\iota(P)=(\cdot,\phi)$ and $P$ enters $e$).
  \label{eq:comb_conf}
  \end{myitem}
The existence of such an $\iota$ follows from~\refeq{bal_paths}.
When $\iota(P)=(\phi,\phi')$, we say that the path $P$ is
{\em attached} to the elements $\phi$ and $\phi'$. So each triangle
in $\Phi$ has three attached paths, by one from each of $H_1,H_2,H_3$,
and each outer edge in $\Phi$ has one attached path.

A converse construction also takes place. More
precisely, consider families $\Phi^+,\Phi^-,\Pscr$ consisting of
copies of some little triangles and edges from $E^+_0(G)$, of copies
of little triangles and edges from $E^-_0(G)$, and of copies
of paths in $H$, respectively. Let $\Phi$ be the disjoint union of
$\Phi^+$ and $\Phi^-$, and $\iota$ a map of $\Pscr$ to
$\Phi^+\times\Phi^-$ satisfying~\refeq{comb_conf}.
We refer to $\Kscr=(\Phi,\Pscr,\iota)$ as a {\em
combinatorial configuration}, or, briefly, a {\em c-configuration}.
Emphasize that we admit some little triangles of $G$ (but not outer
edges) to have copies simultaneously in both $\Phi^+$ and $\Phi^-$.
Now
 \begin{myitem}
define $z(C)$ ($C\in\Cscr(G)$) to be the number of copies of the
triangle $\Delta_C$ in $\Phi^+$ minus the number of copies of
$\Delta_C$ in $\Phi^-$, define $d(e)$ ($e\in E_0(G)$) to be the
number of copies of $e$ in $\Phi^+$ minus the number of copies of $e$
in $\Phi^-$, and define $g:=\sum\{\chi^{P}:P\in\Pscr\}$.
 \label{eq:comb-vect}
  \end{myitem}
Then $z,g,d$ give a v-configuration, denoted by $K(\Kscr)$.
We formally define border $d(\Kscr)$ of $\Kscr$ to be the border of $K
(\Kscr)$. Also we apply to $\Kscr$ adjectives ``proper, essential'' if
$K(\Kscr)$ is such, and similarly for the property of being
``equivalent to'' or ``dominated by'' another configuration.

When no little triangle of $G$ has copies simultaneously in both
$\Phi^+,\Phi^-$, we say that $\Kscr$ is {\em regular}. In particular,
any c-configuration $\Kscr$ representing a v-configuration $K$ by the
first construction is regular; in this case $K(\Kscr)=K$.

\section{Proof of the Theorem} \label{sec:proof}

The proof of Theorem~\ref{tm:puzzle} for a convex grid $G$ falls into
three lemmas. By the arguments in Sections~\ref{sec:lp} and
\ref{sec:flow}, we can deal with c-configurations and, moreover, with
those of them that are regular, essential and proper.

Given a c-configuration $\Kscr=(\Phi,\Pscr,\iota)$,
we say that a little triangle or an outer edge of $G$
or a path of $H$ {\em is in} $\Kscr$ if at least one copy of this
element is contained there. Adding to (deleting from) $\Kscr$
such an element means adding (deleting) exactly one copy of it.
  \Xcomment{
When speaking of a path leaving or entering an element $\phi$ of
$\Phi$ (or an edge in $\phi$), we usually mean the corresponding path of
$\Pscr$ attached to $\phi$ by $\iota$. }

We associate with $\Kscr$ undirected (multi)graph $\Gamma_\Kscr$
whose vertices are the elements of $\Phi$ and whose edges one-to-one
correspond to the paths in $\Pscr$: each path $P\in\Pscr$ generates an
edge connecting $\phi$ and $\phi'$ when $\iota(P)=(\phi,\phi')$
(it is analogous to the graph $\Gamma_\Pi$ associated with a puzzle
$\Pi$, defined in the end of Section~\ref{sec:theorem}).
The (disjoint) union of $\Kscr$ with another or the same
c-configuration $\Kscr'$ is defined in a natural way and denoted by
$\Kscr+\Kscr'$ (its associated graph $\Gamma_{\Kscr+\Kscr'}$ is the
disjoint union of $\Gamma_\Kscr$ and $\Gamma_{\Kscr'}$).

If the interiors of distinct little triangles or edges
$\phi,\phi',\phi''$ of $G$ are traversed by a line of $H$ in this
order, we say that $\phi'$ lies {\em between} $\phi$ and $\phi''$.

We call $\Kscr$ {\em oriented} if all triangles in $\Phi^-$ (the
absorbing triangles) are normal and all triangles in $\Phi^+$ (the
emitting triangles) are turned-over. The first lemma eliminates the
non-oriented configurations.

  \begin{lemma}  \label{lm:orient}
Let a c-configuration $\Kscr=(\Phi,\Pscr,\iota)$ be regular, proper
and essential. There exists a c-configuration $\Kscr'$ such that
either $\Kscr'$ dominates $\Kscr$, or $\Kscr'$ is equivalent to
$\Kscr$ and is oriented.
   \end{lemma}
  \begin{proof}
Since we can consider any regular c-configuration equivalent to
$\Kscr$, one may assume that, among such configurations, $\Kscr$ is
chosen so that
  \begin{myitem}
the number $\eta(\Kscr):=|\Phi|+|\Pscr|$ is as small as possible.
   \label{eq:etaK}
  \end{myitem}

Let us say that a triangle in $\Phi$ is {\em good} if it is either
emitting and turned-over, or absorbing and normal. If all triangles
are good, $\Kscr$ is already oriented. So assume $\Kscr$ contains at
least one bad triangle. Our aim is to show that $\Kscr$ is dominated.

First of all we impose an additional condition on $\Kscr$.
Suppose there is a degenerate path $P\in\Pscr$ attached to a pair of
bad triangles $\Delta\in\Phi^+$ and $\Delta'\in\Phi^-$; so
$\Delta,\Delta'$ share an edge $e$, and $P$ is of the form $v_e$.
Let $e$ be parallel to $\xi_i$ and let $a,a'$ be the
edges of $\Delta,\Delta'$, respectively, parallel to $\xi_{i-1}$.
Observe that $H_{i-1}$ has (one-edge) path $Q$ leaving $v_a$ and
entering $v_{a'}$. When $\Pscr$ contains a copy of $Q$ attached to the
pair $(\Delta,\Delta')$ as well, we call this pair {\em dense}. See
the picture where $i=3$.
 \begin{center}
  \unitlength=1mm
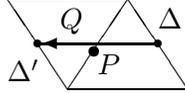
  \begin{picture}(25,12)
  \put(8,0){\line(1,0){16}}
  \put(0,12){\line(1,0){16}}
  \put(8,0){\line(-2,3){8}}
  \put(24,0){\line(-2,3){8}}
  \put(8,0){\line(2,3){8}}
\put(7,8){$Q$}
\put(12,2){$P$}
\put(0,1){$\Delta'$}
\put(20,8){$\Delta$}
\put(4,6){\circle*{1.0}}
\put(11.5,5){\circle*{1.5}}
\put(20,6){\circle*{1.0}}
\thicklines
\put(20,6){\vector(-1,0){15.5}}
  \end{picture}
 \end{center}

We assume that, among all regular c-configurations having the same
border $d(\Kscr)$ and satisfying~\refeq{etaK}, $\Kscr$ is
chosen so that
 \begin{myitem}
the number $\omega(\Kscr)$ of dense pairs in $\Kscr$ is maximum.
  \label{eq:omega}
  \end{myitem}

Suppose the graph $\Gamma_\Kscr$ associated with $\Gamma$ is not
connected. Then $\Kscr$ is the union of two nonempty
c-configurations $\Kscr',\Kscr''$, and we have $d(\Kscr)=d(\Kscr')
+d(\Kscr'')$ and $\eta(\Kscr)=\eta(\Kscr')+\eta(\Kscr')$.
\refeq{etaK} implies that $d(\Kscr')\ne 0$ and $\Kscr'$ is not
equivalent to $\Kscr$. Hence $\Kscr'$ dominates $\Kscr$,
by~\refeq{dominat}(ii). So one may assume that $\Gamma_\Kscr$ is
connected. Then each $\phi\in\Phi$ is reachable in $\Gamma_\Kscr$ by a
path from a vertex representing an outer edge; let $\rho(\phi)$ denote
the minimum number of edges of such a path.

We consider a bad triangle $\Delta$ with $\rho(\Delta)=:\bar\rho$
minimum and proceed by induction on $\bar\rho$. Let $P\in \Pscr$
be a path attached to $\Delta$ and to an element $\phi\in\Phi$
with $\rho(\phi)=\bar\rho-1$. Consider two cases.

{\em Case 1.} Let $\bar\rho=1$. Then $\phi$ is (a copy of) an outer
edge $b$. Assume $b\in E^+_0(G)$; the case $b\in E^-_0(G)$
is symmetric. Then $\Delta$ is absorbing and turned-over, and $P$
leaves $b$ and enters $\Delta$. Let for definiteness $b$ be parallel
to $\xi_2$. For $i=1,2,3$, consider $P_i\in\Pscr$ and
$\phi_i\in\Phi^+$ such that $P_i$ is in $H_i$ and $\iota(P_i)=
(\phi_i,\Delta)$. Let $e_i$ be the edge of $\Delta$
parallel to $\xi_i$. (So $P_2=P$ and $\phi_2=b$.)

Suppose $P_3$ is degenerate, i.e., $P_3=v_{e_3}$. Then $\phi_3$ is a
normal emitting triangle, and therefore, $\phi_3$ is bad.
Take path $Q$ in $H_2$ attached to $\phi_3$, and let
$\iota(Q)=(\phi_3,\tilde\phi)$. The fact that $P_3$ is degenerate
implies that $b,\phi_3,\Delta$ are traversed by a line of $H_2$ in
this order. Hence $H_2$ has path $P'$ leaving $b$ and
entering $\tilde\phi$ and path $Q'$ leaving $\phi_3$ and entering
$\Delta$. Replace in $\Pscr$ the paths $P,Q$ by $P',Q'$,
making $P'$ attached to $b,\tilde\phi$ and making $Q'$ attached to
$\phi_3,\Delta$.

This results in a correct c-configuration $\Kscr'$ with
$\eta(\Kscr')=\eta(K)$ in which $(\phi_3,\Delta)$ becomes a dense
pair. One can see that if $Q$ is nondegenerate, then such a
transformation does not destroy any dense pair of the previous
configuration; so $\omega(\Kscr')>\omega(\Kscr)$,
contradicting~\refeq{omega}. And if $Q$ is degenerate, then
$\phi_3$ and $\tilde\phi$ share an arc of $H_2$, whence $\tilde\phi$
is a turned-over absorbing triangle forming a pair of bad triangles
with $\phi_3$. The only possible dense pair which could be destroyed
by the transformation is just $(\phi_3,\tilde\phi)$
(when this pair is also connected in $\Kscr$ by the corresponding
path in $H_1$). In this case we have $\omega(\Kscr')\ge\omega(\Kscr)$,
so the replacement maintains~\refeq{omega}. Moreover, the new path
leaving $b$ (namely, $P'$) enters a bad triangle (namely, $\tilde\phi$)
as before and is shorter than $P$:
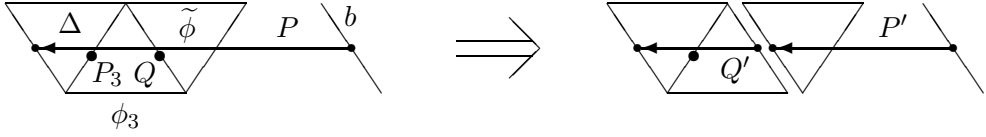
 \begin{center}
  \unitlength=1mm
  \begin{picture}(130,17)
  \put(8,5){\line(1,0){16}}
  \put(0,17){\line(1,0){16}}
  \put(16,17){\line(1,0){16}}
  \put(8,5){\line(-2,3){8}}
  \put(24,5){\line(-2,3){8}}
  \put(8,5){\line(2,3){8}}
  \put(24,5){\line(2,3){8}}
  \put(50,5){\line(-2,3){8}}
\put(7,13){$\Delta$}
\put(23,12.5){$\tilde\phi$}
\put(36,12.5){$P$}
\put(45,14){$b$}
\put(11.5,6.5){$P_3$}
\put(17,6.5){$Q$}
\put(14,1){$\phi_3$}
\put(4,11){\circle*{1.0}}
\put(11.5,10){\circle*{1.5}}
\put(20.5,10){\circle*{1.5}}
\put(46,11){\circle*{1.0}}
\thicklines
\put(46,11){\vector(-1,0){41.5}}
\thinlines
  \put(88,5){\line(1,0){16}}
  \put(80,17){\line(1,0){16}}
  \put(98,17){\line(1,0){16}}
  \put(88,5){\line(-2,3){8}}
  \put(104,5){\line(-2,3){8}}
  \put(106,5){\line(-2,3){8}}
  \put(88,5){\line(2,3){8}}
  \put(106,5){\line(2,3){8}}
  \put(130,5){\line(-2,3){8}}
\put(116,12.5){$P'$}
\put(95,7){$Q'$}
\put(84,11){\circle*{1.0}}
\put(91.5,10){\circle*{1.5}}
\put(100,11){\circle*{1.0}}
\put(102,11){\circle*{1.0}}
\put(126,11){\circle*{1.0}}
\thicklines
\put(126,11){\vector(-1,0){23.5}}
\put(100,11){\vector(-1,0){15.5}}
\thinlines
\put(60,12){\line(1,0){10}}
\put(60,10){\line(1,0){10}}
\put(71,11){\line(-1,1){4}}
\put(71,11){\line(-1,-1){4}}
  \end{picture}
 \end{center}

\noindent Doing so, we eventualy obtain a c-configuration where
$b$ is connected with a bad triangle whose attached path in $H_3$ is
nondegenerate.

Thus, we may assume that $P_3$ is nondegenerate.
Then, by the convexity of $G$, the edge $e_1$ of $\Delta$ does not lie
on the boundary of $G$, and $b$ cannot be the last edge of the side-path
$B^+_2$. We now transform $\Kscr$ as follows. Let $\Delta'$ be the
normal triangle of $G$ containing $e_1$, and $b'$ the edge of
$B^+_2$ next to $b$. Then $H_2$ has path $P'$ leaving $b'$
and entering $\Delta'$ and $H_3$ has path $P'_3$ leaving
$\phi_3$ and entering $\Delta'$ ($P'_3$ is a part of the nondegenerate
$P_3$). We replace in $\Kscr$ the edge $b$ by (one emitting copy of)
$b'$, the triangle $\Delta$ by (one absorbing copy of) $\Delta'$, and
the paths $P,P_3$ by $P',P'_3$, making $P'$ attached to $b',\Delta'$,
and making $P'_3$ attached to $\phi_3,\Delta'$ (while $P_1$ becomes
attached to $\Delta'$ instead of $\Delta$) :
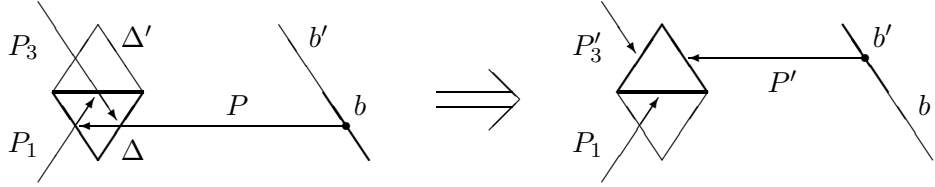
 \begin{center}
  \unitlength=1mm
  \begin{picture}(125,24)
\thicklines
  \put(7,12){\line(1,0){12}}
  \put(13,3){\line(-2,3){6}}
  \put(13,3){\line(2,3){6}}
  \put(49,3){\line(-2,3){6}}
\thinlines
  \put(7,12){\line(2,3){6}}
  \put(19,12){\line(-2,3){6}}
  \put(43,12){\line(-2,3){6}}
   \put(46,7.5){\vector(-1,0){35.5}}
   \put(5,0){\vector(2,3){7.5}}
   \put(5,24){\vector(2,-3){10.5}}
\put(30,9){$P$}
\put(1,4){$P_1$}
\put(1,17){$P_3$}
\put(47,9){$b$}
\put(41,18){$b'$}
\put(16,3){$\Delta$}
\put(16,18){$\Delta'$}
\put(46,7.5){\circle*{1.0}}
%
  \put(88,3){\line(-2,3){6}}
  \put(88,3){\line(2,3){6}}
  \put(124,3){\line(-2,3){6}}
\thicklines
  \put(82,12){\line(1,0){12}}
  \put(82,12){\line(2,3){6}}
  \put(94,12){\line(-2,3){6}}
  \put(118,12){\line(-2,3){6}}
\thinlines
   \put(115,16.5){\vector(-1,0){23.5}}
   \put(80,0){\vector(2,3){7.5}}
   \put(80,24){\vector(2,-3){4.5}}
\put(102,12){$P'$}
\put(76,4){$P_1$}
\put(76,17){$P'_3$}
\put(122,9){$b$}
\put(116,18){$b'$}
\put(115,16.5){\circle*{1.0}}
\put(58,12){\line(1,0){10}}
\put(58,10){\line(1,0){10}}
\put(69,11){\line(-1,1){4}}
\put(69,11){\line(-1,-1){4}}
  \end{picture}
 \end{center}

\noindent This results in a (not necessarily regular) c-configuration
$\Kscr'$ with the border $d(\Kscr)-\chi^b+\chi^{b'}$. By~\refeq{eep},
$\Kscr$ is dominated by $\Kscr'$.

\medskip
{\em Case 2.} Let $\bar\rho>1$. Assume the bad triangle
$\Delta$ is absorbing (and turned-over); the case of emitting
$\Delta$ is symmetric.
Let for definiteness $P$ be in $H_2$, and define
$P_i,\phi_i,e_i$ ($i=1,2,3$) as in Case 1. (So $P=P_2$ and
$\phi=\phi_2$.) Since $\rho(\phi)=\bar\rho-1\ge 1$, $\phi$ is a
good triangle. So $\phi$ is a turned-over emitting triangle and
$P$ is nondegenerate. Arguing as in Case 1, we can impose the
condition that $P_3$ is nondegenerate. This and the convexity
of $G$ imply that neither the edge $e_1$ of $\Delta$
nor the edge $q$ of $\phi$ parallel to $\xi_1$ is in $b(G)$. Let
$\Delta'$ be the normal little triangle of $G$ containing $e_1$, and
$\phi'$ the normal triangle containing $q$.
We replace $\Delta,\phi$ in $\Phi$ by $\Delta',\phi'$.

More precisely, when $\Delta$ is replaced by $\Delta'$, we
accordingly replace the paths $P,P_3$ attached to $\Delta$ by
paths $P',P'_3$ (while $P_1$ preserves, becoming attached to
$\Delta'$). Here $P'$ is the path of $H_2$ leaving
$\phi'$ and entering $\Delta'$, and $P'_3$ is the path of $H_3$
leaving $\phi_3$ and entering $\Delta'$ (as before, $P'_3$ is a part
of the nondegenerate path $P_3$). And when replacing $\phi$ by
$\phi'$, we should also replace path $\tilde Q$ of $H_3$ attached to
$\phi$, entering triangle $\tilde\phi\in\Phi^-$ say, by path
$\tilde Q'$ of $H_3$ leaving $\phi'$ and entering $\tilde\phi$.
($\tilde Q'$ exists since $\phi$ lies between $\phi'$ and
$\tilde\phi$.) The path of $H_1$ attached to
$\phi$ becomes attached to $\phi'$. This gives a c-configuration
$\Kscr'$ in which the added triangle $\phi'$ is bad and has the rank
$\rho(\phi')$ equal to $\bar\rho-1$.

We have $d(\Kscr')=d(\Kscr)$ and $\eta(\Kscr')=\eta(\Kscr)$.
The latter implies that $\Kscr'$ is regular, i.e., $\Kscr$ has no
emitting copy of $\Delta'$ or $\phi'$. For otherwise, cancelling in
$\Kscr'$ one emitting copy and one absorbing copy of the same little
triangle of $G$ and properly concatenating their attached paths, we
would obtain a configuration with a smaller value of $\eta$, contrary
to~\refeq{etaK}. Finally, one can see that neither $\Delta$ nor
$\phi$ can be involved in dense pairs of $\Kscr$. Hence
no dense pair is destroyed while constructing $\Kscr'$, implying
$\omega(\Kscr')\ge \omega(\Kscr)$. Now the result follows by induction
on $\bar\rho$.
  \end{proof}

Thus, it suffices to consider only oriented configurations.

A puzzle $\Pi=(\Fscr,\Pscr)$ generates an oriented
c-configuration $(\Phi,\Pscr,\iota)$ in a natural way: $\Phi^+$ is
the set of turned-over triangles in $\Fscr$ and edges in $b^+(\Pi)$,
$\Phi^-$ is the set of normal triangles in $\Fscr$ and edges in
$b^-(\Pi)$, and for each $u-v$~path $P\in\Pscr$, $\iota(P)$ is
the pair $(\phi\in\Phi^+,\phi\in\Phi^-)$ such that the point $u$ is
contained in $\phi$ and $v$ is contained in $\phi'$. Such a
{\em puzzle c-configuration} is denoted by $\Kscr_\Pi$.

The next lemma describes a situation when an oriented configuration
can be split into two configurations one of which is a puzzle
configuration. Let us say that paths $P,P'$ of $H$
are {\em crossing} if they are not parallel and their interiors
have a point in common, and that $P$ and a little triangle $\Delta$ of $G$
are {\em overlapping} if $P$ meets the interior of $\Delta$:
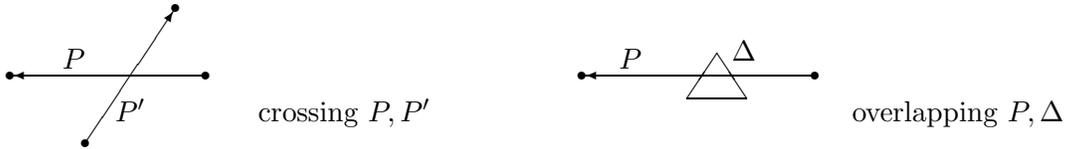
 \begin{center}
  \unitlength=1mm
  \begin{picture}(140,18)
  \put(26,9){\vector(-1,0){25.5}}
  \put(10,0){\vector(2,3){11.7}}
    \put(0,9){\circle*{1.0}}
    \put(26,9){\circle*{1.0}}
    \put(10,0){\circle*{1.0}}
    \put(22,18){\circle*{1.0}}
    \put(7,10){$P$}
    \put(14,3){$P'$}
    \put(33,3){crossing $P,P'$}
  \put(107,9){\vector(-1,0){30.5}}
  \put(90,6){\line(1,0){8}}
  \put(98,6){\line(-2,3){4}}
  \put(90,6){\line(2,3){4}}
    \put(76,9){\circle*{1.0}}
    \put(107,9){\circle*{1.0}}
    \put(81,10){$P$}
    \put(96,11){$\Delta$}
    \put(112,3){overlapping $P,\Delta$}
  \end{picture}
 \end{center}

One can see that the puzzle configurations are precisely those having
neither crossing nor overlapping pairs.
Given an oriented c-configuration $\Kscr=(\Phi,\Pscr,\iota)$,
define its {\em minimal pre-configuration} $\Kmin=
(\Psi,\Pmin,\widehat\iota)$ as follows. Let $\Psi^+$ (resp. $\Psi^-$)
be the set of little triangles and outer edges of $G$ having at
least one copy in $\Phi^+$ (resp. $\Phi^-$). Then
$\Psi=\Psi^+\cup\Psi^-$. The set $\Pmin$ is formed by taking for each
edge $e\in E(G)$ contained in a member of $\Psi^+$, one
(inclusion-wise) minimal path in $\Pscr$ with the beginning $v_e$,
taking for each edge $e\in E(G)$ contained in a member of $\Psi^-$,
one minimal path in $\Pscr$ with the end $v_e$, and
ignoring repeated paths if arise. Define $\widehat\iota$ to be the map
attaching a $u-v$~path $P\in\Pmin$ to the pair $(\phi\in\Psi^+,
\phi'\in\Psi^-)$ such that $u\in\phi$ and $v\in\phi'$ (this pair
is unique since $K$ is oriented). Note that $\Kmin$ need not be a
c-configuration since some triangles (outer edges) in it may have more
than three (resp. one) attached paths.
  \begin{lemma}  \label{lm:min-conf}
Let a c-configuration $\Kscr=(\Phi,\Pscr,\iota)$ be proper and
oriented, and let $\Kmin=(\Psi,\Pmin,\widehat\iota)$ be its minimal
pre-configuration.
Suppose $\Kmin$ contains neither crossing paths nor overlapping a path
and a triangle. Then: (a) $\Kmin$ is a puzzle c-configuration, and
(b) $\Kmin$ either is equivalent to $\Kscr$ or dominates $\Kscr$.
   \end{lemma}
   \begin{proof}
From the non-existence of paths in $\Pmin$ overlapping triangles in
$\Psi$ it easily follows that for each element $\phi\in \Psi^+$
and each edge $e$ in $\phi$, there is exactly one path $P\in \Pmin$
leaving $e$, and similarly for each element $\phi'\in \Psi^-$
and each edge $e'$ in $\phi'$, there is exactly one path $P'\in \Pmin$
entering $e'$. Hence $\Kmin$ is a c-configuration, and now the absence
of crossing paths in $\Kscr$ implies that it is a puzzle
configuration, yielding (a). Next, one can rearrange the attaching
map $\iota$ in $\Kscr$ so that $\Kscr$ be represented as the union of
$\Kmin$ and some c-configuration $\Kscr''$. This implies (b),
by~\refeq{dominat}(ii).
   \end{proof}

For $i=1,2,3$, a sequence $(\phi_1,\ldots,\phi_k)$ of distinct
little triangles or edges of $G$ is called an $i$-{\em chain} if
their interiors are traversed in this order by a path of $H_i$.
If $(\Delta,\Delta')$ is an $i$-chain of two normal little triangles
and there is no normal triangle between them, we say that $\Delta$
is the $i$-{\em predecessor} of $\Delta'$, and similarly for
turned-over triangles.

Our final lemma is the following.
   \begin{lemma}  \label{lm:cross-overlap}
Let a c-configuration $\Kscr=(\Phi,\Pscr,\iota)$ be proper, essential
and oriented. If $\Kscr$ is not equivalent to a puzzle
c-configuration, then $\Kscr$ is dominated.
   \end{lemma}
   \begin{proof}
Since we can replace $\Kscr$ by any oriented c-configuration
equivalent to $\Kscr$ (e.g., by taking the union of $r$ copies of
$\Kscr$ for any $r$), one may assume that, among such configurations,
$\Kscr$ is chosen so that:
  \begin{myitem1}
\item[(i)]
there are sufficiently many copies of each member of $\Phi\cup\Pscr$;
  \item[(ii)] subject to (i), the number $t(\Kscr)$ of little
triangles of $G$ having copies in $\Phi$ is maximum;
  \item[(iii)] subject to (i),(ii), the number $p(\Kscr)$ of paths of
$H$ having copies in $\Kscr$ is maximum.
  \label{eq:eta2}
  \end{myitem1}

From (iii) it follows that
  \begin{myitem}
for any (not necessarily distinct) vertices $u_1,u_2,u_3,u_4$
occurring in a path of $H$ in this order, if $\Pscr$ contains
copies of both $u_1-u_3$ path $P$ and $u_2-u_4$ path $P'$,
then $\Pscr$ contains copies of both $u_1-u_4$ path $Q$ and
$u_2-u_3$ path $Q'$ as well, and vice versa.
  \label{eq:u1-u4}
  \end{myitem}
Indeed, if at least one of $Q,Q'$ is not in $\Pscr$, we
can add $Q,Q'$ to $\Pscr$ and delete $P,P'$ from $\Pscr$, accordingly
correcting the map $\iota$. This increases $p(\Kscr)$. (Recall that
adding to $\Kscr$ a triangle or an outer edge of $G$ or a path of $H$
means adding {\em one} copy of this element, and similarly for
deleting an element.) The reverse assertion is proved similarly.

\smallskip
Also we assume that the minimal pre-configuration $\Kmin$ contains
crossing paths or overlapping a path and a triangle; otherwise the
result immediately follows from Lemma~\ref{lm:min-conf}. We show that
$\Kscr$ is dominated in both cases.

\medskip
{\em Case 1}. Let $\Kmin$ contain crossing a $u-v$ path $P$ and a
$u'-v'$ path $Q$. Assume for definiteness that $P$ is in $H_2$ and
minimal among the paths of $\Pscr$ beginning at $u$, and that $Q$
is in $H_1$ ($P$ is anti-parallel to $\xi_1$ and $Q$ is anti-parallel
to $\xi_3$); the case when $P$ is minimal among the paths ending
at $v$ is symmetric. Observe that the point $w$ where $P,Q$ intersect
is a vertex of $H_2$.
Let $\Delta$ be the normal little triangle whose edge parallel to
$\xi_2$ contains $w$ as the median point. Then $\Delta$ is not in
$\Phi$. For otherwise $\Pscr$ would contain a path from some
vertex $w'$ to $w$ (as $\Delta$ is absorbing). Applying~\refeq{u1-u4}
to $w',u,w,v$ or to $u,w',w,v$, we obtain that $\Pscr$ contains the
$u-w$ path, contradicting the minimality of $P$.

Next we proceed as follows. For $i=1,2,3$, let $e_i$ be the edge
of $\Delta$ parallel to $\xi_i$. (So $w=v_{e_2}$.)
Take the turned-over little triangle $\nabla$ containing $e_3$.
Let $e'_1,e'_2$ be the edges of $\nabla$ parallel to $\xi_1,\xi_2$,
respectively. Then $H_2$ has $u-w$ path $P'$ and
$v_{e'_2}-v$ path $P''$, and $H_1$ has $u'-v_{e_1}$ path
$Q'$:
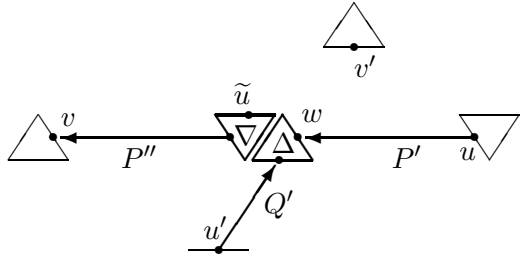
 \begin{center}
  \unitlength=1mm
  \begin{picture}(70,33)
  \put(0,12){\begin{picture}(8,6)%
     \put(0,0){\line(1,0){8}}         
     \put(0,0){\line(2,3){4}}
     \put(8,0){\line(-2,3){4}}
            \end{picture}}
  \put(42,27){\begin{picture}(8,6)%
     \put(0,0){\line(1,0){8}}         
     \put(0,0){\line(2,3){4}}
     \put(8,0){\line(-2,3){4}}
            \end{picture}}
  \put(24,0){\line(1,0){8}}         
  \put(60,12){\begin{picture}(8,6)%
     \put(0,6){\line(1,0){8}}         
     \put(4,0){\line(2,3){4}}
     \put(4,0){\line(-2,3){4}}
            \end{picture}}
  \put(32.5,12){\begin{picture}(8,6)%
     \thicklines
     \put(0,0){\line(1,0){8}}         
     \put(0,0){\line(2,3){4}}
     \put(8,0){\line(-2,3){4}}
     \thinlines
            \end{picture}}
  \put(27.5,12){\begin{picture}(8,6)%
     \thicklines
     \put(0,6){\line(1,0){8}}         
     \put(4,0){\line(2,3){4}}
     \put(4,0){\line(-2,3){4}}
     \thinlines
            \end{picture}}
    \put(6,15){\circle*{1.0}}
    \put(29.5,15){\circle*{1.0}}
    \put(38.5,15){\circle*{1.0}}
    \put(62,15){\circle*{1.0}}
    \put(28,0){\circle*{1.0}}
    \put(36,12){\circle*{1.0}}
    \put(32,18){\circle*{1.0}}
    \put(46,27){\circle*{1.0}}
     \thicklines
  \put(29.5,15){\vector(-1,0){22.5}}
  \put(62,15){\vector(-1,0){22.5}}
  \put(28,0){\vector(2,3){7.5}}
    \thinlines
    \put(7,16.5){$v$}
    \put(15,11){$P''$}
    \put(30,14){$\nabla$}
    \put(35,13){$\Delta$}
    \put(39,17){$w$}
    \put(51,11){$P'$}
    \put(60,12){$u$}
    \put(26,1.5){$u'$}
    \put(34,5){$Q'$}
    \put(46,23){$v'$}
    \put(30,19.5){$\tilde u$}
  \end{picture}
 \end{center}

\medskip
Add (one copy of) the triangle $\Delta$ to $\Phi^-$, the triangle
$\nabla$ to $\Phi^+$, and the paths $P',P'',Q'$ together with the
degenerate path $v_{e_3}$ (``connecting'' $\Delta$ and $\nabla$ in
$H_3$) to $\Pscr$. Accordingly delete $P,Q$ from $\Pscr$.
The attachments for the added elements are assigned in a natural way
(e.g., $\iota(P'):=(\phi,\Delta)$, where $\phi$ is the element of the
old $\Phi^+$ to which $P$ was attached). This increases the value of
the parameter $t$ (since $\Delta$ is added while the new $\Kscr$
contains a copy of each triangles from the previous $\Kscr$,
by assumption~\refeq{eta2}(i)). However, $\Kscr$ becomes an
``incomplete'' configuration since $\nabla$ has no attached path in
$H_1$, and similarly for element $\widehat\phi$ of $\Phi^-$ to which
$Q$ was attached. We cannot improve $\Kscr$ straightforwardly because
the points $\tilde u:=v_{e'_1}$ and $v'$ do not lie on one line of
$H_1$.

Our aim is to improve this $\Kscr$ without decreasing the current
value of $t$, in order to obtain a correct c-configuration $\Kscr'$
either dominating or equivalent to the initial $\Kscr$. This will
yield the result in the former case and lead to a contradiction with
assumption~\refeq{eta2}(ii) in the latter case.

First of all we iteratively construct a sequence $S$ of alternating
members of $\Phi$ and $\Pscr$ as follows. Start with $\Delta_1:=
\widehat\phi$. Let $\Delta_i\in\Phi$ be the last element of the current
$S$. If $\Delta_i$ is an outer edge, halt.
Otherwise add $P_{i+1},\Delta_{i+1}$ to $S$,
where $P_{i+1}$ is attached to $\Delta_i,\Delta_{i+1}$.
More precisely: (a) if $i$ is odd (and $\Delta_i$ is a
normal triangle), then $P_{i+1}$ is a path of $H_2$ and
$\iota(P_{i+1})=(\Delta_{i+1},\Delta_i)$, and (b) if $i$ is even (and
$\Delta_i$ is a turned-over triangle), then $P_{i+1}$ is a path of
$H_1$ and $\iota(P_{i+1})=(\Delta_i,\Delta_{i+1})$. Let $\Delta_{q+1}$
be the last element of the final $S$. Clearly the edge $b:=
\Delta_{q+1}$ belongs to $B^+_2$ when $q$ is odd, and to $B^-_1$ when
$q$ is even.

Assume $q$ is odd; the case of $q$ even is examined analogously.
For $i=1,\ldots,q$, let $Q_i\in\Pscr$ be the path of $H_3$ attached
to $\Delta_i$ (it enters $\Delta_i$ for $i$ odd, and leaves
$\Delta_i$ for $i$ even). Let $\Delta'_i$ be the other element of
$\Phi$ to which $Q_i$ is attached. We say that the triangle
$\Delta_i$ is {\em squeezed} if $i$ is odd and $Q_i$ is degenerate.

We first explain how to transform $\Kscr$ into the desired correct
c-configuration when no $\Delta_i$ is squeezed.
By the convexity of $G$ (and regardless of the squeezedness of any
$\Delta_i$), the line in the plane parallel to $\xi_3$
and passing the point $\tilde u$ separates $S$ from $B^+_3$ (letting
$B^+_3$ be the common vertex of $B^-_2$ and $B^-_1$ when they meet).
This implies that $S$ can be shifted by distance 1 in
the direction of $\xi_2$ (approaching $B^+_3$).
More precisely, each triangle $\Delta_i$ has 3-predecessor $\tilde
\Delta_i$ in $G$, and $B^+_2$ contains edge $\tilde b$ next to $b$.
See the picture where $q=3$.
 \begin{center}
  \unitlength=1mm
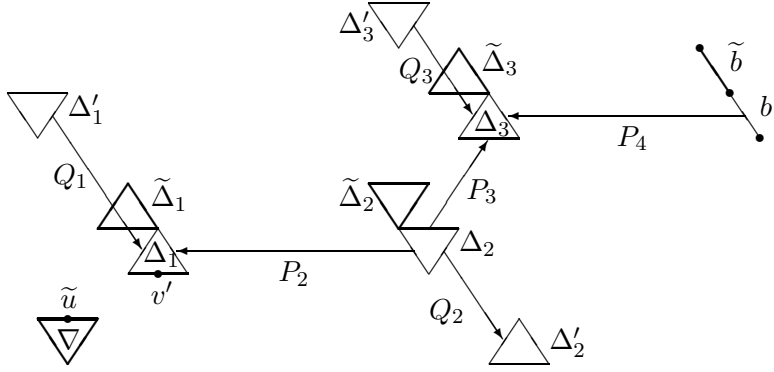
  \begin{picture}(110,48)
  \put(64,0){\begin{picture}(8,6)%
     \put(0,0){\line(1,0){8}}         
     \put(0,0){\line(2,3){4}}
     \put(8,0){\line(-2,3){4}}
            \end{picture}}
  \put(16,12){\begin{picture}(8,6)%
     \put(0,0){\line(1,0){8}}         
     \put(0,0){\line(2,3){4}}
     \put(8,0){\line(-2,3){4}}
            \end{picture}}
  \put(60,30){\begin{picture}(8,6)%
     \put(0,0){\line(1,0){8}}         
     \put(0,0){\line(2,3){4}}
     \put(8,0){\line(-2,3){4}}
            \end{picture}}
  \put(0,30){\begin{picture}(8,6)%
     \put(0,6){\line(1,0){8}}         
     \put(4,0){\line(2,3){4}}
     \put(4,0){\line(-2,3){4}}
            \end{picture}}
  \put(52,12){\begin{picture}(8,6)%
     \put(0,6){\line(1,0){8}}         
     \put(4,0){\line(2,3){4}}
     \put(4,0){\line(-2,3){4}}
            \end{picture}}
  \put(48,42){\begin{picture}(8,6)%
     \put(0,6){\line(1,0){8}}         
     \put(4,0){\line(2,3){4}}
     \put(4,0){\line(-2,3){4}}
            \end{picture}}
  \put(12,18){\begin{picture}(8,6)%
     \thicklines
     \put(0,0){\line(1,0){8}}         
     \put(0,0){\line(2,3){4}}
     \put(8,0){\line(-2,3){4}}
     \thinlines
            \end{picture}}
  \put(56,36){\begin{picture}(8,6)%
     \thicklines
     \put(0,0){\line(1,0){8}}         
     \put(0,0){\line(2,3){4}}
     \put(8,0){\line(-2,3){4}}
     \thinlines
            \end{picture}}
  \put(4,0){\begin{picture}(8,6)%
     \thicklines
     \put(0,6){\line(1,0){8}}         
     \put(4,0){\line(2,3){4}}
     \put(4,0){\line(-2,3){4}}
     \thinlines
            \end{picture}}
  \put(48,18){\begin{picture}(8,6)%
     \thicklines
     \put(0,6){\line(1,0){8}}         
     \put(4,0){\line(2,3){4}}
     \put(4,0){\line(-2,3){4}}
     \thinlines
            \end{picture}}
     \put(100,30){\line(-2,3){4}}
     \thicklines
     \put(96,36){\line(-2,3){4}}
     \thinlines
    \put(100,30){\circle*{1.0}}
    \put(96,36){\circle*{1.0}}
    \put(92,42){\circle*{1.0}}
    \put(8,6){\circle*{1.0}}
    \put(20,12){\circle*{1.0}}
  \put(54,15){\vector(-1,0){31.5}}
  \put(98,33){\vector(-1,0){31.5}}
  \put(6,33){\vector(2,-3){11.8}}
  \put(58,15){\vector(2,-3){7.8}}
  \put(54,45){\vector(2,-3){7.8}}
  \put(56,18){\vector(2,3){7.8}}
    \put(6.5,2){$\nabla$}
    \put(18,13.5){$\Delta_1$}
    \put(60,15){$\Delta_2$}
    \put(62,31){$\Delta_3$}
    \put(8,33){$\Delta'_1$}
    \put(72,2){$\Delta'_2$}
    \put(44,44){$\Delta'_3$}
    \put(19,21){$\tilde\Delta_1$}
    \put(44,21){$\tilde\Delta_2$}
    \put(63,39){$\tilde\Delta_3$}
    \put(100,33){$b$}
    \put(96,39){$\tilde b$}
    \put(36,11){$P_2$}
    \put(61,22){$P_3$}
    \put(81,29){$P_4$}
    \put(6,24){$Q_1$}
    \put(56,6){$Q_2$}
    \put(52,38){$Q_3$}
    \put(7,7.5){$\tilde u$}
    \put(19,8){$v'$}
  \end{picture}
 \end{center}

These triangles $\tilde\Delta_i$ and the elements
$\tilde\Delta_0:=\nabla$ and
$\tilde\Delta_{q+1}:=\tilde b$ are connected in $H$ by paths
$P'_1,\ldots,P'_{q+1}$ in a natural way: $P'_i$ is the path of $H_1$
leaving $\tilde\Delta_i$ and entering $\tilde\Delta_{i+1}$ when
$i$ is odd, and the path of $H_2$ leaving $\tilde\Delta_{i+1}$ and
entering $\tilde\Delta_i$ when $i$ is even.
Also there are paths $Q'_1,\ldots,Q'_q$ of $H_3$ such that $Q'_i$
leaves $\Delta'_i$ and enters $\tilde\Delta_i$ when $i$ is odd
(as $\Delta_i$ is not squeezed, and therefore, $\tilde\Delta_i$ lies
between $\Delta'_i$ and $\Delta_i$), and
$Q'_i$ leaves $\tilde\Delta_i$ and enters $\Delta'_i$ when $i$ is
even.

Add to $\Kscr$ the triangles $\tilde\Delta_1,\ldots,\tilde\Delta_q$,
the paths $P'_1,\ldots,P'_{q+1},Q'_1,\ldots,Q'_q$ and the edge
$\tilde b$, making $P'_i$ attached to $\tilde\Delta_{i-1},
\tilde\Delta_i$, and making $Q'_j$ attached to
$\tilde\Delta_j,\Delta'_j$. Accordingly delete from $\Kscr$ the
triangles $\Delta_1,\ldots,\Delta_q$, the paths $P_2,\ldots,P_{q+1},
Q_1,\ldots,Q_q$ and the outer edge $b$. This results in a correct
c-configuration $\Kscr'$. Moreover, $\Kscr'$ has the border $d(\Kscr)-
\chi^b+\chi^{\tilde b}$. Therefore, $\Kscr'$ dominates the initial
$\Kscr$, by~\refeq{eep}.

\smallskip
Next suppose there is a squeezed $\Delta_i$ ($i$ is odd); let $i$ be
minimum among such triangles. Form the triangles $\tilde\Delta_0,
\ldots,\tilde\Delta_{i-1}$ and paths
$P'_1,\ldots,P'_i,Q'_1,\ldots,Q'_{i-1}$ as above.
Take paths $R,D\in\Pscr$ attached to $\Delta'_i$ and belonging
to $H_2$ and $H_1$, respectively. Let $\phi,\phi'$ be the other
(normal) triangles to which $R,D$ are attached, respectively.
Since $\Delta_i$ is squeezed, $(\Delta_{i+1},\Delta_i,\Delta'_i,\phi)$
is a 2-chain and $(\tilde\Delta_{i-1},\Delta'_i,\phi')$ is a 1-chain:
 \begin{center}
  \unitlength=1mm
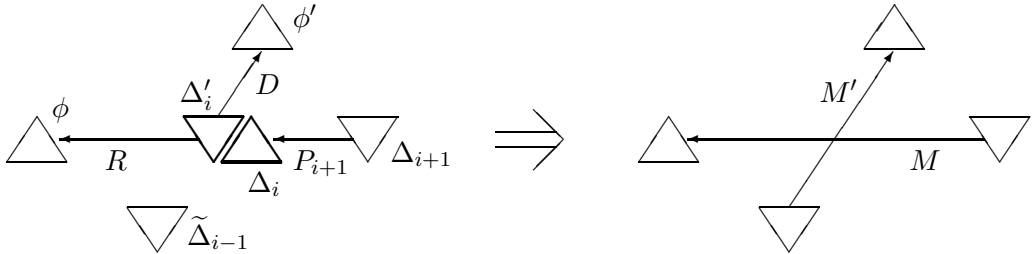
  \begin{picture}(140,33)
  \put(0,12){\begin{picture}(8,6)%
     \put(0,0){\line(1,0){8}}         
     \put(0,0){\line(2,3){4}}
     \put(8,0){\line(-2,3){4}}
            \end{picture}}
  \put(30,27){\begin{picture}(8,6)%
     \put(0,0){\line(1,0){8}}         
     \put(0,0){\line(2,3){4}}
     \put(8,0){\line(-2,3){4}}
            \end{picture}}
  \put(16,0){\begin{picture}(8,6)%
     \put(0,6){\line(1,0){8}}         
     \put(4,0){\line(2,3){4}}
     \put(4,0){\line(-2,3){4}}
            \end{picture}}
  \put(44,12){\begin{picture}(8,6)%
     \put(0,6){\line(1,0){8}}         
     \put(4,0){\line(2,3){4}}
     \put(4,0){\line(-2,3){4}}
            \end{picture}}
  \put(28.5,12){\begin{picture}(8,6)%
     \thicklines
     \put(0,0){\line(1,0){8}}         
     \put(0,0){\line(2,3){4}}
     \put(8,0){\line(-2,3){4}}
     \thinlines
            \end{picture}}
  \put(23.5,12){\begin{picture}(8,6)%
     \thicklines
     \put(0,6){\line(1,0){8}}         
     \put(4,0){\line(2,3){4}}
     \put(4,0){\line(-2,3){4}}
     \thinlines
            \end{picture}}
  \put(25.5,15){\vector(-1,0){18.5}}
  \put(46,15){\vector(-1,0){10.5}}
  \put(28,18){\vector(2,3){5.8}}
    \put(6,18){$\phi$}
    \put(13,11){$R$}
    \put(23,20){$\Delta'_i$}
    \put(38,11){$P_{i+1}$}
    \put(51,12){$\Delta_{i+1}$}
    \put(24,1){$\tilde\Delta_{i-1}$}
    \put(32,8){$\Delta_i$}
    \put(33,21){$D$}
    \put(38,30){$\phi'$}
  \put(84,12){\begin{picture}(8,6)%
     \put(0,0){\line(1,0){8}}         
     \put(0,0){\line(2,3){4}}
     \put(8,0){\line(-2,3){4}}
            \end{picture}}
  \put(114,27){\begin{picture}(8,6)%
     \put(0,0){\line(1,0){8}}         
     \put(0,0){\line(2,3){4}}
     \put(8,0){\line(-2,3){4}}
            \end{picture}}
  \put(100,0){\begin{picture}(8,6)%
     \put(0,6){\line(1,0){8}}         
     \put(4,0){\line(2,3){4}}
     \put(4,0){\line(-2,3){4}}
            \end{picture}}
  \put(128,12){\begin{picture}(8,6)%
     \put(0,6){\line(1,0){8}}         
     \put(4,0){\line(2,3){4}}
     \put(4,0){\line(-2,3){4}}
            \end{picture}}
  \put(130,15){\vector(-1,0){39.5}}
  \put(104,6){\vector(2,3){13.8}}
    \put(120,11){$M$}
    \put(108,20){$M'$}
     \put(65,14){\line(1,0){8}}         
     \put(65,16){\line(1,0){8}}
     \put(74,15){\line(-1,1){4}}
     \put(74,15){\line(-1,-1){4}}
  \end{picture}
 \end{center}

\smallskip
Let $M$ be the path of $H_2$ leaving $\Delta_{i+1}$ and entering
$\phi$, and $M'$ the path of $H_1$ leaving $\tilde\Delta_{i-1}$
and entering $\phi'$. We add to $\Kscr$ the triangles
$\tilde\Delta_1,\ldots,\tilde\Delta_{i-1}$
and the paths $P'_1,\ldots,P'_i,Q'_1,\ldots, Q'_{i-1},M,M'$
and accordingly delete the triangles $\Delta_1,\ldots,
\Delta_i$ and $\Delta'_i$ and the paths $P_2,\ldots,P_{i+1},Q_1,
\ldots,Q_i,R,D$.
(Note that if $P_i$ is degenerate, then $\tilde\Delta_{i-1}$ and
$\Delta'_i$ are copies of the same triangle of $G$; we consider them
as different objects one of which is added and the other is deleted.)
The resulting $\Kscr'$ is a correct c-configuration with the same
border $d(\Kscr)$. But $t(\Kscr')>t(\Kscr)$ (as $\Delta$ was added,
while deleting the above triangles does not affect $t$,
by~\refeq{eta2}(i)). This contradicts~\refeq{eta2}(ii).

\medskip
{\em Case 2.}
Let $\Kmin$ contain overlapping a path $P$ and a triangle $\phi$.
One may assume that $P$ is a $u-v$~path of $H_1$ and that
$P$ is minimal among the paths in $\Pscr$ beginning at $u$.
Let $\iota(P)=(\phi',\phi'')$. Then $\phi$ lies between $\phi'$ and
$\phi''$. Notice that there is no normal (absorbing) triangle
$\tilde\phi\in\Phi$ between $\phi'$ and $\phi''$.
For if such a $\tilde\phi$ exists, then the end vertex $w$ of the path
of $H_1$ attached to $\tilde\phi$ is an intermediate vertex of $P$.
But then $\Pscr$ contains the $u-w$~path (by~\refeq{u1-u4}), contrary
to the minimality of $P$. So $\phi$ is a turned-over (emitting)
triangle.

Take path $Q$ of $H_2$ attached to $\phi$; let
$\iota(Q)=(\phi,\psi)$. Since the absorbing element $\psi$ cannot be a
normal triangle lying between $\phi'$ and $\phi''$
(by the argument above), the path $Q$ is nondegenerate.
Let $e$ be the edge of $\phi$ parallel to $\xi_2$, and $\Delta$
the normal little triangle of $G$ containing $e$.
Then $\Delta$ lies between $\phi'$ and $\phi''$; let $P'$ be the path
of $H_1$ leaving $\phi'$ and entering $\Delta$. Note that $\Kscr$
contains no copy of $\Delta$ (again by the argument above).
Next, let $e'$ be the edge of $\Delta$ parallel to $\xi_3$,
and $\nabla$ the turned-over triangle of $G$ containing $e'$. Then
$\nabla$ lies between $\phi$ and $\psi$ (as $Q$ is
nondegenerate); let $Q'$ be the path of $H_2$ leaving $\nabla$ and
entering $\psi$.
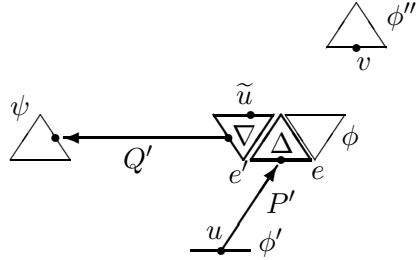
 \begin{center}
  \unitlength=1mm
  \begin{picture}(50,33)
  \put(0,12){\begin{picture}(8,6)%
     \put(0,0){\line(1,0){8}}         
     \put(0,0){\line(2,3){4}}
     \put(8,0){\line(-2,3){4}}
            \end{picture}}
  \put(42,27){\begin{picture}(8,6)%
     \put(0,0){\line(1,0){8}}         
     \put(0,0){\line(2,3){4}}
     \put(8,0){\line(-2,3){4}}
            \end{picture}}
  \put(24,0){\line(1,0){8}}         
  \put(36.5,12){\begin{picture}(8,6)%
     \put(0,6){\line(1,0){8}}         
     \put(4,0){\line(2,3){4}}
     \put(4,0){\line(-2,3){4}}
            \end{picture}}
  \put(32,12){\begin{picture}(8,6)%
     \thicklines
     \put(0,0){\line(1,0){8}}         
     \put(0,0){\line(2,3){4}}
     \put(8,0){\line(-2,3){4}}
     \thinlines
            \end{picture}}
  \put(27,12){\begin{picture}(8,6)%
     \thicklines
     \put(0,6){\line(1,0){8}}         
     \put(4,0){\line(2,3){4}}
     \put(4,0){\line(-2,3){4}}
     \thinlines
            \end{picture}}
    \put(6,15){\circle*{1.0}}
    \put(29,15){\circle*{1.0}}
    \put(28,0){\circle*{1.0}}
    \put(36,12){\circle*{1.0}}
    \put(32,18){\circle*{1.0}}
    \put(46,27){\circle*{1.0}}
     \thicklines
  \put(29,15){\vector(-1,0){22}}
  \put(28,0){\vector(2,3){7.5}}
    \thinlines
    \put(26,1.5){$u$}
    \put(46,24){$v$}
    \put(29.5,14){$\nabla$}
    \put(34.5,13){$\Delta$}
    \put(44,14){$\phi$}
    \put(33,0){$\phi'$}
    \put(50,30){$\phi''$}
    \put(0,18){$\psi$}
    \put(34,5){$P'$}
    \put(15,11){$Q'$}
    \put(40,9.5){$e$}
    \put(29,9){$e'$}
    \put(30,19.5){$\tilde u$}
  \end{picture}
 \end{center}

\medskip
Add to $\Kscr$ the triangles $\Delta,\nabla$, the paths $P',Q'$,
the degenerate path $v_e$ (``connecting'' $\phi$ and $\Delta$ in
$H_2$) and the degenerate path $v_{e'}$ (``connecting'' $\nabla$ and
$\Delta$ in $H_3$), assigning the attachments for them in an obvious
way. Accordingly delete from $\Kscr$ the paths $P,Q$. This
results in an ``incomplete'' c-configuration, but having
a larger value of $t$, in which $\nabla$ and $\phi''$ have no attached
paths in $H_1$. (It cannot be improved straightforwardly since
$v$ and the median point $\tilde u$ of the edge of $\nabla$ parallel to
$\xi_1$ do not lie on one line of $H_1$).
So we have a situation as in Case 1 and proceed in a similar
way to transform $\Kscr$ into a correct c-configuration $\Kscr'$
either dominating the initial $\Kscr$ or being equivalent to $\Kscr$
but having a larger value of $t$.

This completes the proof of the lemma.
\end{proof}

By Lemmas~\ref{lm:orient} and~\ref{lm:cross-overlap}, any
non-dominated proper essential configuration is equivalent to a puzzle
configuration. This implies Theorem~\ref{tm:puzzle}, in view of
explanations in Sections~\ref{sec:lp},\ref{sec:flow}.

\bigskip
\noindent{\bf Remark.} Analysing the proof of
Lemma~\ref{lm:cross-overlap}, one sees that, in fact, a slightly
sharper version of this lemma is obtained; namely, taking into account
assumption~\refeq{eta2}(ii) and the construction of the minimal
pre-configuration $\Kmin$:
  \begin{myitem}
if a c-configuration $\Kscr$ is proper, essential and oriented and if
$\Kscr$ is not dominated, then $\Kscr$ is equivalent to a puzzle
configuration $\Kscr_\Pi$ such that the set of triangles of the puzzle
$\Pi$ includes all little triangles of $G$ having copies in $\Kscr$.
  \label{eq:sharp}
  \end{myitem}

\section{Concluding Remarks} \label{sec:concl}

We conclude this paper with several remarks.

\smallskip
{\em First},
for a cocirculation $h$ in $G$ and a tandem $\tau=(e,e')$, call
$\delta_h(\tau):=h(e)-h(e')$ the {\em discrepancy} of $h$ at $\tau$.
So $h$ is concave if the discrepancy at each tandem is nonnegative.
A more general problem ($\ast$) is to find a cocirculation
$h$ having a given border $\sigma$ and obeying prescribed lower bounds
$c$ on the discrepancies: $\delta_h(\tau)\ge c(\tau)$ for each
$\tau\in\Tscr(G)$. This is reduced to the case of
zero bounds when $c$ comes up from another cocirculation $g$ in $G$.
More precisely, let $c(\tau):=\delta_g(\tau)$ for each tandem $\tau$.
Re-define the required border by $\sigma'(e):=\sigma(e)-g(e)$
for each outer edge $e$. Then $h'$ is a concave cocirculation with
the border $\sigma'$ if and only if $h:=h'+g$ is a cocirculation with
the border $\sigma$ satisfying the lower bound $c$ on the
discrepancies. Thus, the corresponding changes in the puzzle
inequalities~\refeq{gall-ineq} and in the monotone
condition~\refeq{monot} give a solvability criterion for
problem~($\ast$) with a cocirculation-induced $c$.

In particular, the puzzle criterion modified in this way works when
all tandem discrepancies are required to be greater than or equal to a
prescribed constant $\alpha\in\Rset$. This is because there exists a
cocirculation $g$ in $G$ where the discrepancy at each tandem is
exactly $\alpha$. (Such a $g$ is constructed easily: assuming w.l.o.g.
that $G$ is a 3-side grid of size $n$, put $g(e_i):=(k-2i+1)\alpha$
($i=1,\ldots,k$) for each side-path $(e_1,\ldots,e_k)$ in $G$.)

\smallskip
{\em Second},
from the sharper version of Lemma~\ref{lm:cross-overlap} given
in~\refeq{sharp} one derives that each puzzle $\Pi$ determining a
facet of $\BC(G)$ is (uniquely) determined by its boundary $b(\Pi)$.

Indeed, suppose $\Pi_1,\Pi_2$ are two different puzzles with
$b(\Pi_1)=b(\Pi_2)$. Let $\Kscr_i$ stand for the c-configuration
induced by $\Pi_i$; one may assume that $\Kscr_i$ is proper and
essential. Then $\Kscr:=\Kscr_1+\Kscr_2$ is an oriented
c-configuration equivalent to $\Kscr_i$.
Assume $\Kscr$ is not dominated and take the puzzle $\Pi$ as
in~\refeq{sharp}. We have $b(\Kscr_\Pi)=b(\Kscr_i)$, so the number
$q$ of triangles in $\Pi$ is equal to the number $q_1$ of triangles in
$\Pi_1$, by~\refeq{par_puzzle}.
On the other hand, the fact that $\Pi_1$ and $\Pi_2$ are different
implies that $\Kscr$ involves more little triangles of $G$ compared
with $\Kscr_1$. This implies $q>q_1$, by~\refeq{sharp}; a
contradiction.

\begin{figure}[tb]
 \begin{center}
  \unitlength=1mm
  \begin{picture}(155,54)
      \put(0,0){\begin{picture}(80,54)
  \put(0,0){\line(1,0){72}}
  \put(0,0){\line(2,3){36}}
  \put(72,0){\line(-2,3){36}}
  \put(48,0){\begin{picture}(8,6)%
     \thicklines
     \put(0,0){\line(1,0){8}}         
     \put(0,0){\line(2,3){4}}
     \put(8,0){\line(-2,3){4}}
            \end{picture}}
  \put(28,6){\begin{picture}(8,6)%
     \thicklines
     \put(0,0){\line(1,0){8}}         
     \put(0,0){\line(2,3){4}}
     \put(8,0){\line(-2,3){4}}
            \end{picture}}
  \put(8,12){\begin{picture}(8,6)%
     \thicklines
     \put(0,0){\line(1,0){8}}         
     \put(0,0){\line(2,3){4}}
     \put(8,0){\line(-2,3){4}}
            \end{picture}}
  \put(24,24){\begin{picture}(8,6)%
     \thicklines
     \put(0,0){\line(1,0){8}}         
     \put(0,0){\line(2,3){4}}
     \put(8,0){\line(-2,3){4}}
            \end{picture}}
  \put(44,18){\begin{picture}(8,6)%
     \thicklines
     \put(0,0){\line(1,0){8}}         
     \put(0,0){\line(2,3){4}}
     \put(8,0){\line(-2,3){4}}
            \end{picture}}
  \put(40,36){\begin{picture}(8,6)%
     \thicklines
     \put(0,0){\line(1,0){8}}         
     \put(0,0){\line(2,3){4}}
     \put(8,0){\line(-2,3){4}}
            \end{picture}}
  \put(40,6){\begin{picture}(8,6)%
     \thicklines
     \put(0,6){\line(1,0){8}}         
     \put(4,0){\line(2,3){4}}
     \put(4,0){\line(-2,3){4}}
            \end{picture}}
  \put(20,12){\begin{picture}(8,6)%
     \thicklines
     \put(0,6){\line(1,0){8}}         
     \put(4,0){\line(2,3){4}}
     \put(4,0){\line(-2,3){4}}
            \end{picture}}
  \put(36,24){\begin{picture}(8,6)%
     \thicklines
     \put(0,6){\line(1,0){8}}         
     \put(4,0){\line(2,3){4}}
     \put(4,0){\line(-2,3){4}}
            \end{picture}}
    \put(4,0){\circle*{1.0}}
    \put(28,0){\circle*{1.0}}
    \put(52,0){\circle*{1.0}}
    \put(50,3){\circle*{1.0}}
    \put(70,3){\circle*{1.0}}
    \put(32,6){\circle*{1.0}}
    \put(30,9){\circle*{1.0}}
    \put(34,9){\circle*{1.0}}
    \put(42,9){\circle*{1.0}}
    \put(46,9){\circle*{1.0}}
    \put(12,12){\circle*{1.0}}
    \put(44,12){\circle*{1.0}}
    \put(10,15){\circle*{1.0}}
    \put(14,15){\circle*{1.0}}
    \put(22,15){\circle*{1.0}}
    \put(26,15){\circle*{1.0}}
    \put(24,18){\circle*{1.0}}
    \put(48,18){\circle*{1.0}}
    \put(46,21){\circle*{1.0}}
    \put(50,21){\circle*{1.0}}
    \put(58,21){\circle*{1.0}}
    \put(28,24){\circle*{1.0}}
    \put(26,27){\circle*{1.0}}
    \put(30,27){\circle*{1.0}}
    \put(38,27){\circle*{1.0}}
    \put(42,27){\circle*{1.0}}
    \put(40,30){\circle*{1.0}}
    \put(22,33){\circle*{1.0}}
    \put(44,36){\circle*{1.0}}
    \put(42,39){\circle*{1.0}}
    \put(46,39){\circle*{1.0}}
    \put(34,51){\circle*{1.0}}
     \thicklines
  \put(4,0){\vector(2,3){7.7}}
  \put(28,0){\vector(2,3){3.7}}
  \put(44,12){\vector(2,3){3.7}}
  \put(24,18){\vector(2,3){3.7}}
  \put(40,30){\vector(2,3){3.7}}
  \put(70,3){\vector(-1,0){15.5}}
  \put(42,9){\vector(-1,0){7.5}}
  \put(22,15){\vector(-1,0){7.5}}
  \put(58,21){\vector(-1,0){7.5}}
  \put(38,27){\vector(-1,0){7.5}}
  \put(46,9){\vector(2,-3){3.7}}
  \put(26,15){\vector(2,-3){3.7}}
  \put(42,27){\vector(2,-3){3.7}}
  \put(22,33){\vector(2,-3){3.7}}
  \put(34,51){\vector(2,-3){7.7}}
       \end{picture}}
      \put(80,0){\begin{picture}(80,54)
  \put(0,0){\line(1,0){72}}
  \put(0,0){\line(2,3){36}}
  \put(72,0){\line(-2,3){36}}
  \put(48,0){\begin{picture}(8,6)%
     \thicklines
     \put(0,0){\line(1,0){8}}         
     \put(0,0){\line(2,3){4}}
     \put(8,0){\line(-2,3){4}}
            \end{picture}}
  \put(32,12){\begin{picture}(8,6)%
     \thicklines
     \put(0,0){\line(1,0){8}}         
     \put(0,0){\line(2,3){4}}
     \put(8,0){\line(-2,3){4}}
            \end{picture}}
  \put(8,12){\begin{picture}(8,6)%
     \thicklines
     \put(0,0){\line(1,0){8}}         
     \put(0,0){\line(2,3){4}}
     \put(8,0){\line(-2,3){4}}
            \end{picture}}
  \put(28,18){\begin{picture}(8,6)%
     \thicklines
     \put(0,0){\line(1,0){8}}         
     \put(0,0){\line(2,3){4}}
     \put(8,0){\line(-2,3){4}}
            \end{picture}}
  \put(36,18){\begin{picture}(8,6)%
     \thicklines
     \put(0,0){\line(1,0){8}}         
     \put(0,0){\line(2,3){4}}
     \put(8,0){\line(-2,3){4}}
            \end{picture}}
  \put(40,36){\begin{picture}(8,6)%
     \thicklines
     \put(0,0){\line(1,0){8}}         
     \put(0,0){\line(2,3){4}}
     \put(8,0){\line(-2,3){4}}
            \end{picture}}
  \put(36,12){\begin{picture}(8,6)%
     \thicklines
     \put(0,6){\line(1,0){8}}         
     \put(4,0){\line(2,3){4}}
     \put(4,0){\line(-2,3){4}}
            \end{picture}}
  \put(28,12){\begin{picture}(8,6)%
     \thicklines
     \put(0,6){\line(1,0){8}}         
     \put(4,0){\line(2,3){4}}
     \put(4,0){\line(-2,3){4}}
            \end{picture}}
  \put(32,18){\begin{picture}(8,6)%
     \thicklines
     \put(0,6){\line(1,0){8}}         
     \put(4,0){\line(2,3){4}}
     \put(4,0){\line(-2,3){4}}
            \end{picture}}
    \put(4,0){\circle*{1.0}}
    \put(28,0){\circle*{1.0}}
    \put(52,0){\circle*{1.0}}
    \put(50,3){\circle*{1.0}}
    \put(70,3){\circle*{1.0}}
    \put(36,12){\circle*{1.0}}
    \put(34,15){\circle*{1.0}}
    \put(38,15){\circle*{1.0}}
    \put(42,15){\circle*{1.0}}
    \put(12,12){\circle*{1.0}}
    \put(10,15){\circle*{1.0}}
    \put(14,15){\circle*{1.0}}
    \put(30,15){\circle*{1.0}}
    \put(32,18){\circle*{1.0}}
    \put(40,18){\circle*{1.0}}
    \put(42,21){\circle*{1.0}}
    \put(58,21){\circle*{1.0}}
    \put(30,21){\circle*{1.0}}
    \put(34,21){\circle*{1.0}}
    \put(38,21){\circle*{1.0}}
    \put(36,24){\circle*{1.0}}
    \put(22,33){\circle*{1.0}}
    \put(44,36){\circle*{1.0}}
    \put(42,39){\circle*{1.0}}
    \put(46,39){\circle*{1.0}}
    \put(34,51){\circle*{1.0}}
     \thicklines
  \put(4,0){\vector(2,3){7.7}}
  \put(28,0){\vector(2,3){7.7}}
  \put(36,24){\vector(2,3){7.7}}
  \put(70,3){\vector(-1,0){15.5}}
  \put(30,15){\vector(-1,0){15.5}}
  \put(58,21){\vector(-1,0){15.5}}
  \put(42,15){\vector(2,-3){7.7}}
  \put(22,33){\vector(2,-3){7.7}}
  \put(34,51){\vector(2,-3){7.7}}
       \end{picture}}
  \end{picture}
 \end{center}
\caption{Two puzzles with equal boundaries} \label{fig:two_gall}
  \end{figure}
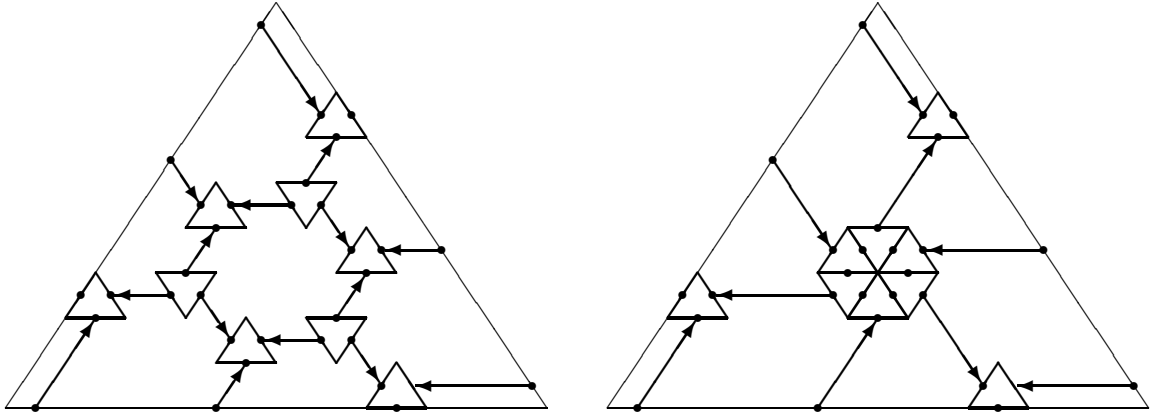

\smallskip
Different puzzles with equal boundaries do exist. An example for a
3-side grid is shown in Fig.~\ref{fig:two_gall}.

A puzzle determined by its boundary is called {\em rigid}. Knutson,
Tao and Woodward proved that in the case of 3-side grids the
facet-determining puzzles are exactly the rigid ones.
They, moreover, obtained a combinatorial characterization of
facet-determining puzzles, implying that they are recognizable in
polynomial time.
\begin{theorem} {\rm\cite{KTW}}  \label{tm:facet}
Let $\Pi=(\Fscr,\Pscr)$ be a puzzle in a 3-side grid $G$ such that
$\Fscr$ is nonempty and different from the set of all little triangles
of $G$. The following are equivalent:
  \begin{itemize}
\item[(i)] $\Pi$ determines a facet of $\BC(G)$;
\item[(ii)] $\Pi$ is rigid;
\item[(iii)] $\Pi$ admits no gentle circuits.
  \end{itemize}
  \end{theorem}

To explain the notion of gentle path/circuit, let $R$ be the set of
little rhombi of $G$ that are split by a path in $\Pscr$ into two
parallelograms. Let $G_0$ be the subgraph of $G$ induced by
the edges separating either a triangle in $\Fscr$ and a rhombus in
$R$ (the {\em tp-edges}), or a rhombus in $R$ and a little triangle
contained in no member of $\Fscr\cup R$ (the {\em pn-edges}).
Re-orient each tp-edge (resp. pn-edge) $e$ so that the triangle of
$\Fscr$ (resp. the rhombus in $R$) containing $e$ lie on the right.
A path or circuit $P$ of $G_0$ is called {\em gentle} if, when moving
along $P$ from an edge to the next edge, the angle of turn is either
$0^\circ$ or $60^\circ$, never $120^\circ$. For example, the circuit
surrounding the hexagon formed by the six central triangles in the
right puzzle in Fig.~\ref{fig:two_gall} is gentle.

One can show that Theorem~\ref{tm:facet} remains valid for an arbitrary
convex grid $G$. (Implication (i)$\to$(ii) has already been shown. The
method of proof of (ii)$\to$(iii) and (iii)$\to$(i) given in~\cite{KTW}
is applicable to an arbitrary convex grid, as it, in essense, does not
depend on the shape of the convex region $\Rscr$ spanned by $G$.
Roughly speaking, the proof of (ii)$\to$(iii) relies on a local
transformation of a puzzle $\Pi$ having a gentle circuit $C$.
It creates another puzzle with the same boundary by re-arranging $\Pi$
only within the 1-neighbourhood of $C$ (being the union of little
triangles and rhombi sharing common edges with $C$).
The proof of (iii)$\to$(i)
uses the function on the tp- and pn-edges whose value
on an edge $e$ is the number of all maximal gentle paths with the
first edge $e$. When $\Pi$ has no gentle circuits, this function
(regardless of the shape of $\Rscr$) is easily transformed into a
concave cocirculation $h_0$ in $G$ for which the tandem inequality
is strict on each little rhombus separated by a tp- or pn-edge.
Using $h_0$, it is routine to construct $|E_0(G)|-2$
concave cocirculations whose borders are linearly independent and
orthogonal to the border of $\Kscr_\Pi$.) We omit details of the
proof here.

It is not difficult to check that any puzzle $\Pi$ with
$\Fscr=\emptyset$ and $|\Pscr|=1$ is facet-determining as well (such a
puzzle can arise when $\Rscr$ has $\ge 4$ sides). When $\Fscr=\emptyset$
and $|\Pscr|\ge 2$, $\Pi$ is already not facet-determining as it is
the union of two disjoint puzzles.

\smallskip
{\em Third}, by a result of Knutson and Tao~\cite{KT} on integral
honeycombs, a feasible integer-valued function $\sigma$ on the outer
edges of a 3-side grid $G$ is extendable to an integer concave
cocirculation. In a forthcoming paper~\cite{Kar} we show that this is
generalized to a convex grid $G$ and, moreover, a sharper property
takes place: a concave cocirculation $h$ in a convex grid $G$ can be
turned into an integer concave cocirculation preserving the values of
$h$ on all outer edges $e$ with $h(e)\in\Zset$ and on the edges of all
little triangles where $h$ is integral for each of the three edges.

\medskip
{\bf Acknowledgement.} I am thankful to Vladimir Danilov and Gleb
Koshevoy for many stimulating discussions on discrete concave
functions and related topics.

\end{document}